\documentclass[reqno,12pt]{amsart}

\newtheorem{Theorem}{Theorem}[section]

\newtheorem{Lemma}[Theorem]{Lemma}

\newtheorem{Conjecture}[Theorem]{Conjecture}
\theoremstyle{remark}
\newtheorem{Remark}[Theorem]{Remark}
\numberwithin{equation}{section}

\allowdisplaybreaks

\title[Noncommutative hypergeometric series]
{Summation formulae for\\noncommutative hypergeometric series}

\author[Michael Schlosser]{Michael Schlosser$^*$}
\thanks{$^*$Supported by Austrian Science Fund FWF, grant P17563-N13,
and partly by EC's IHRP Programme, grant HPRN-CT-2001-00272
``Algebraic Combinatorics in Europe''.}
\address{Fakult\"at f\"ur Mathematik, Universit\"at Wien,
Nordbergstra{\ss}e 15, A-1090 Wien, Austria}
\email{schlosse@ap.univie.ac.at}
\urladdr{http://www.mat.univie.ac.at/{\textasciitilde}schlosse}

\date{November 6, 2004}

\begin{document}

\begin{abstract}
We establish several summation formulae for hypergeometric and basic
hypergeometric series involving noncommutative parameters and argument.
These results were inspired by a recent paper of
J.~A.~Tirao [Proc.\ Nat.\ Acad.\ Sci.~100 (14) (2003), 8138--8141].
\end{abstract}

\maketitle

\section{Introduction}

Hypergeometric series with noncommutative parameters and argument,
in the special case involving square matrices,
have recently been studied by a number of researchers including
(in alphabetical order) Dur\'an, Duval, Gr\"unbaum, Iliev, Ovsienko,
Pacharoni, Tirao, and others. See \cite{DG,DO,G,GI,GPT1,GPT2,GPT3,T}
for some selected papers. The subject of hypergeometric series
involving matrices is closely related to and partly overlapping the
theory of orthogonal matrix polynomials. The study of the latter was
initiated by Krein~\cite{Kr} and subsequently has experienced a steady
development. Whereas a good amount of theory of orthogonal matrix
polynomials has already been worked out, see e.g.\ Dur\'an and
Van Assche~\cite{DV}, Dur\'an and L\'opez-Rodr\'{\i}guez~\cite{DL},
and Tirao~\cite{T}, it seems as appropriate to study
noncommutative hypergeometric series (involving not only matrices
but more generally arbitrary noncommutative parameters of some unit ring,
or, in the case of infinite series, of some Banach algebra) from an
entirely elementary point of view.
This includes the search for identities for noncommutative
hypergeometric and noncommutative basic hypergeometric series,
extending their classical commutative versions which
can be found, for instance, in the standard textbooks of
Bailey~\cite{B}, Slater~\cite{S}, Gasper and Rahman~\cite{GR},
and of Andrews, Askey, and Roy~\cite{AAR}.

This paper contains some results of our search which we hope
will be the starting point of a systematic study towards a theory
of identities for noncommutative hypergeometric series and
their basic analogues ($q$-analogues). The special types of
noncommutative hypergeometric series we are considering were inspired
by a recent paper of Tirao~\cite{T}. To be precise, we consider
noncommutative hypergeometric series of types I and II.
Tirao's matrix extension of the Gau{\ss} hypergeometric function
belongs to type I, according to our terminology in
Section~\ref{secprel}. We would like to stress that by
``noncommutative'' we do not mean ``$q$-commutative''
or ``quasi commutative'', i.e.\ involving a relation like $yx=qxy$.
For some results on the latter, see the papers by Koornwinder~\cite{Ko}
and Volkov~\cite{V} and the references therein.
Unless we specify explicit commutation
relations (which will sometimes happen)
our parameters, elements of an abstract noncommutative unit ring,
are understood {\em not} to commute with each other.

Our paper is organized as follows. In Section~\ref{secprel}
we first recall the classical definitions for
($q$-)shifted factorials and (basic) hypergeometric series,
and then define their noncommutative versions.
In Section~\ref{secelid} we prove by induction a couple of
lemmas containing simple addition formulae for the noncommutative
($Q$-)shifted factorials. Sections~\ref{sec2f1}, \ref{sec3f2} and
\ref{sec7f6} are the heart pieces of our paper. Here we derive
noncommutative extensions of several important terminating summations for
hypergeometric and basic hypergeometric series, in particular of the
($q$-)Chu--Vandermonde summations, the ($q$-)Pfaff--Saalsch\"utz
summation, and of Dougall's ${}_7F_6$ summation.
Concerning the latter summation, we were unfortunately not able to
establish a noncommutative $Q$-analogue. With other words, the problem
of finding a noncommutative $Q$-Dougall
summation (or Jackson summation) is still open.
The summations in Sections~\ref{sec2f1}--\ref{sec7f6}
are proved by entirely elementary means, namely by induction.
The situation is quite different in Section~\ref{secnt} where,
now working in some abstract Banach algebra,
we give some nonterminating identities,
in particular, two $Q$-Gau{\ss} summations (derived using a formal
argument) as conjectures, and further two nonterminating $Q$-binomial
theorems which we prove using functional equations.
Finally, in Section~\ref{secmore} we indicate two ways for obtaining
even more identities for noncommutative (basic) hypergeometric series,
one of them is ``telescoping'', the other is ``reversing all products''.

The results is this paper provide (to the best of our knowledge) a
first collection of identities for noncommutative hypergeometric and
noncommutative basic hypergeometric series. We were honestly surprised
when finding the identities in this paper which extend some of the most
important summation formulae in the theory of (basic) hypergeometric
series to noncommuting parameters. A continuation of our program may
include the derivation of yet other summations but also noncommutative
analogues of some of the classical transformation formulae
(as listed in \cite[Appendix~III]{B} and \cite[Appendix~III]{GR}).
Another issue left open is the study of the ``type II''
Gau{\ss} hypergeometric function and the ``type I'' and ``type II''
$Q$-Gau{\ss} hypergeometric functions from the view-point of the
second-order ($Q$-)differential equations they (presumably) satisfy
when considered in an appropriate analytic setting, in the spirit of
Tirao's~\cite{T} illuminating investigation of the (what we call)
ordinary (i.e.\ ``non-$Q$'') ``type I'' case.
While this paper offers an elementary approach (the terminating
summations are obtained by two applications of induction in each
instance) and focuses on explicit summation formulae only,
we feel that in order to gain more insight the subject
of noncommutative ($Q$-)hypergeometric series should also be
investigated from a broader perspective, connecting it to other
theories (in combinatorics, representation theory, physics, etc.)
where possible and appropriate.
We have strong confidence that our summation formulae will be useful
in the theory of ($Q$-)orthogonal matrix polynomials (e.g., to show
that certain ``principal'' specializations of these polynomials factor)
and may even provide motivation for defining new selected
families of ($Q$-)orthogonal matrix polynomials. Furthermore,
it seems not too farfetched to expect that our identities will have
applications in some noncommutative models of mathematical physics
(see also \cite{DO}).

The author would like to thank George Gasper and Hjalmar Rosengren for
their comments. We are especially indebted to the latter for providing
an explicit matrix solution to a particular set of algebraic equations
(see Remark~\ref{rem}).

\section{Preliminaries}\label{secprel}

\subsection{Classical (commutative) hypergeometric and
basic hypergeometric series}

The standard references for hypergeometric series
and basic hypergeometric series are \cite{S} and \cite{GR}, respectively.

Define the {\em shifted factorial} for all integers $k$ by
the following quotient of gamma functions,
\begin{equation*}
(a)_k:=\frac{\Gamma(a+k)}{\Gamma(a)}.
\end{equation*}
Further, the (ordinary) {\em hypergeometric $_{r+1}F_r$ series}
is defined as
\begin{equation}\label{defohyp}
{}_{r+1}F_r\!\left[\begin{matrix}a_1,a_2,\dots,a_{r+1}\\
b_1,b_2,\dots,b_r\end{matrix};z\right]:=
\sum _{k=0}^{\infty}\frac {(a_1)_k\dots(a_{r+1})_k}
{(b_1)_k\dots(b_r)_k}\frac{z^k}{k!}.
\end{equation}

Let $q$ be a complex number such that $0<|q|<1$. Define the
{\em $q$-shifted factorial} for all integers $k$ (including infinity) by
\begin{equation*}
(a;q)_k:=\prod_{j=1}^k(1-aq^j).
\end{equation*}
We write
\begin{equation}\label{defhyp}
{}_{r+1}\phi_r\!\left[\begin{matrix}a_1,a_2,\dots,a_{r+1}\\
b_1,b_2,\dots,b_r\end{matrix}\,;q,z\right]:=
\sum _{k=0}^{\infty}\frac {(a_1;q)_k\dots (a_{r+1};q)_k}
{(b_1;q)_k\dots(b_r;q)_k}\frac{z^k}{(q;q)_k},
\end{equation}
to denote the {\em basic hypergeometric ${}_{r+1}\phi_r$ series}.
In \eqref{defohyp} and \eqref{defhyp}, $a_1,\dots,a_{r+1}$ are
called the {\em upper parameters}, $b_1,\dots,b_r$ the
{\em lower parameters}, $z$ is the {\em argument}, and
(in \eqref{defhyp}) $q$ the {\em base} of the series.
The ${}_{r+1}\phi_r$ series in \eqref{defhyp} reduces to
the ${}_{r+1}F_r$ series in \eqref{defohyp} after first replacing
all parameters $a_i$ by $q^{a_i}$ and  $b_i$ by $q^{b_i}$ and
then letting $q\to 1$. This possibility of taking limits to obtain
ordinary hypergeometric series from basic hypergeometric series
is not shared by the noncommutative versions of \eqref{defohyp}
and \eqref{defhyp} which we will define in Subsection~\ref{subsecnc}.
   
The hypergeometric ${}_{r+1}F_r$ series terminates if one of
the upper parameters, say $a_{r+1}$,
is of the form $-n$, for a nonnegative integer $n$.
On the other hand, the basic hypergeometric ${}_{r+1}\phi_r$ series
terminates if one of the upper parameters, say $a_{r+1}$, is of the
form $q^{-n}$, for a nonnegative integer $n$.
See \cite[p.~45]{S} and \cite[p.~25]{GR} for the criteria of
when the hypergeometric, respectively basic hypergeometric, series
converge if they do not terminate.

The classical theories of hypergeometric basic hypergeometric series
contain several important summation and transformation formulae
involving ${}_{r+1}F_r$ and ${}_{r+1}\phi_r$ series.
Many of these summation theorems require
that the parameters satisfy the condition of being
either balanced and/or very-well-poised.
An ${}_{r+1}F_r$ hypergeometric series is called
{\em balanced} (or {\em $1$-balanced}) if
$b_1+\cdots+b_r=a_1+\cdots+a_{r+1}+1$ and $z=1$.
More generally, it is called  {\em $k$-balanced} if
$b_1+\cdots+b_r=a_1+\cdots+a_{r+1}+k$ and $z=1$.
Similarly, an ${}_{r+1}\phi_r$ basic hypergeometric series is
called {\em balanced} if $b_1\cdots b_r=a_1\cdots a_{r+1}q$ and $z=q$.
An ${}_{r+1}F_r$ series is called {\em well-poised} if
$a_1+1=a_2+b_1=\cdots=a_{r+1}+b_r$ and is {\em very-well-poised}
if in addition $a_2=\frac{a_1}2+1$.
Note that this choice of $a_2$ entails that the factor
\begin{equation*}
\frac{\frac{a_1}2+k}{\frac{a_1}2}
\end{equation*}
appears in a very-well-poised series.
Similarly, an ${}_{r+1}\phi_r$ basic hypergeometric series is
called {\em well-poised} if $a_1q=a_2b_1=\cdots=a_{r+1}b_r$ and is
{\em very-well-poised} if in addition $a_2=-a_3=q\sqrt{a_1}$.
Here this choice of $a_2$ and $a_3$ entails that the factor
\begin{equation*}
\frac {1-a_1q^{2k}}{1-a_1}
\end{equation*}
appears in a very-well-poised basic series.
In both cases (ordinary and basic), the parameter
$a_1$ is referred to as the
{\em special parameter} of the very-well-poised series.

\subsection{Noncommutativity}\label{subsecnc}

Let $R$ be a unit ring (i.e., a ring with a multiplicative identity).
Throughout this article, the elements of $R$ shall be denoted by
capital letters $A,B,C,\dots$. In general these elements do not commute
with each other; however, we may sometimes specify certain commutation
relations explicitly. We denote the identity by $I$ and the zero element
by $O$.
Whenever a multiplicative inverse element exists for any $A\in R$,
we denote it by $A^{-1}$. (Since $R$ is a unit ring, we have
$AA^{-1}=A^{-1}A=I$.) On the other hand, as we shall implicitly assume
that all the expressions which appear are well defined, whenever
we write $A^{-1}$ we assume its existence.
For instance, in \eqref{defncpochI} and \eqref{defncpochII}
we assume that $C_i+jI$ is invertible
for all $1\le i\le r$, $0\le j<k$.

An important special case is when $R$ is the ring of $n\times n$
square matrices (our notation is certainly suggestive with respect
to this interpretation), or, more generally, one may view
$R$ as a space of some abstract operators.

Let $\mathbb Z$ be the set of integers.
For $l,m\in\mathbb Z\cup\{\pm\infty\}$ we define
the noncommutative product as follows:
\begin{equation}\label{ncprod}
\prod_{j=l}^mA_j=\begin{cases}
1&m=l-1\\
A_lA_{l+1}\dots A_m&m\ge l\\
A_{l-1}^{-1}A_{l-2}^{-1}\dots A_{m+1}^{-1}&m<l-1
\end{cases}.
\end{equation}
Note that
\begin{equation}\label{invprod}
\prod_{j=l}^mA_j=\prod_{j=m+1}^{l-1}A_{m+l-j}^{-1},
\end{equation}
for all $l,m\in\mathbb Z\cup\{\pm\infty\}$.
We will make use of \eqref{invprod} at the end of Section~\ref{sec2f1}
when reversing the order of summation of a series and pulling out factors.

Let $k\in\mathbb Z$. We define the generalized
{\em noncommutative shifted factorial of type I}\/ by
\begin{equation}\label{defncpochI}
\left\lceil\begin{matrix}
A_1,A_2,\dots,A_r\\C_1,C_2,\dots,C_r
\end{matrix};Z\right\rfloor_k:=
\prod_{j=1}^k\left[\left(\prod_{i=1}^r
(C_i+(k-j)I)^{-1}(A_i+(k-j)I)\right)Z\right],
\end{equation}
and the {\em noncommutative shifted factorial of type II}\/ by
\begin{equation}\label{defncpochII}
\left\lfloor\begin{matrix}
A_1,A_2,\dots,A_r\\C_1,C_2,\dots,C_r
\end{matrix};Z\right\rceil_k:=
\prod_{j=1}^k\left[\left(\prod_{i=1}^r
(C_i+(j-1)I)^{-1}(A_i+(j-1)I)\right)Z\right].
\end{equation}
Note the unusual usage of brackets (``floors'' and ``ceilings''
are intermixed) on the left-hand sides of \eqref{defncpochI} and
\eqref{defncpochII} which is intended to suggest that the products
involve noncommuting factors in a prescribed order.
In both cases, the product,
read from left to right, starts with a denominator factor.
The brackets in the form ``$\lceil-\rfloor$'' are intended to denote
that the factors are {\em falling},
while in ``$\lfloor-\rceil$'' that they are {\em rising}.

If $Z=I$, we write
\begin{equation}\label{defncpochIa}
\left\lceil\begin{matrix}
A_1,A_2,\dots,A_r\\C_1,C_2,\dots,C_r
\end{matrix};I\right\rfloor_k=\left\lceil\begin{matrix}
A_1,A_2,\dots,A_r\\C_1,C_2,\dots,C_r
\end{matrix}\right\rfloor_k,
\end{equation}
and
\begin{equation}\label{defncpochIIa}
\left\lfloor\begin{matrix}
A_1,A_2,\dots,A_r\\C_1,C_2,\dots,C_r
\end{matrix};I\right\rceil_k=\left\lfloor\begin{matrix}
A_1,A_2,\dots,A_r\\C_1,C_2,\dots,C_r
\end{matrix}\right\rceil_k,
\end{equation}
for simplicity in notation.

We define the {\em noncommutative hypergeometric series of type I}\/ by
\begin{equation}\label{defnchypI}
{}_{r+1}F_r\!\left\lceil\begin{matrix}
A_1,A_2,\dots,A_{r+1}\\C_1,C_2,\dots,C_r
\end{matrix};Z\right\rfloor:=
\sum_{k\ge 0}
\left\lceil\begin{matrix}
A_1,A_2,\dots,A_{r+1}\\C_1,C_2,\dots,C_r,I
\end{matrix}\,;Z\right\rfloor_k,
\end{equation}
and the
{\em noncommutative hypergeometric series of type II}\/ by
\begin{equation}\label{defnchypII}
{}_{r+1}F_r\!\left\lfloor\begin{matrix}
A_1,A_2,\dots,A_{r+1}\\C_1,C_2,\dots,C_r
\end{matrix};Z\right\rceil:=
\sum_{k\ge 0}
\left\lfloor\begin{matrix}
A_1,A_2,\dots,A_{r+1}\\C_1,C_2,\dots,C_r,I
\end{matrix}\,;Z\right\rceil_k.
\end{equation}
In each case, the series terminates if one of the upper
parameters $A_i$ is of the form $-nI$. The situation is more delicate
if the series is nonterminating. In this case we shall assume
that $R$ is a Banach algebra with norm $\|\cdot\|$.
Then the series converges in $R$ if $\|Z\|<1$.
If $\|Z\|=1$ the series may converge in $R$ for some particular
choice of upper and lower parameters. Exact conditions depend
on the Banach algebra $R$.

Throughout this paper, $Q$ will be a parameter which commutes with
any of the other parameters appearing in the series.
(For instance, a central element such as $Q=qI$, a scalar multiple
of the unit element in $R$,
for $qI\in R$, trivially satisfies this requirement.)

Let $k\in\mathbb Z$.
The generalized {\em noncommutative $Q$-shifted factorial of type I}\/
is defined by
\begin{equation}\label{defncpochIQ}
\left\lceil\begin{matrix}
A_1,A_2,\dots,A_r\\C_1,C_2,\dots,C_r
\end{matrix};Q,Z\right\rfloor_k:=
\prod_{j=1}^k\left[\left(\prod_{i=1}^r
(I-C_iQ^{k-j})^{-1}(I-A_iQ^{k-j})\right)Z\right].
\end{equation}
Similarly, the generalized
{\em noncommutative $Q$-shifted factorial of type II}\/ is defined by
\begin{equation}\label{defncpochIIQ}
\left\lfloor\begin{matrix}
A_1,A_2,\dots,A_r\\C_1,C_2,\dots,C_r
\end{matrix};Q,Z\right\rceil_k:=
\prod_{j=1}^k\left[\left(\prod_{i=1}^r
(I-C_iQ^{j-1})^{-1}(I-A_iQ^{j-1})\right)Z\right].
\end{equation}
Formally, one may let $k\to\infty$
(or, if one desires, even $k\to-\infty$).

We define the {\em noncommutative basic
hypergeometric series of type I}\/ by
\begin{equation}\label{defnchypIQ}
{}_{r+1}\phi_r\!\left\lceil\begin{matrix}
A_1,A_2,\dots,A_{r+1}\\C_1,C_2,\dots,C_r
\end{matrix};Q,Z\right\rfloor:=
\sum_{k\ge 0}
\left\lceil\begin{matrix}
A_1,A_2,\dots,A_{r+1}\\C_1,C_2,\dots,C_r,Q
\end{matrix};Q,Z\right\rfloor_k,
\end{equation}
and the {\em noncommutative basic
hypergeometric series of type II}\/ by
\begin{equation}\label{defnchypIIQ}
{}_{r+1}\phi_r\!\left\lfloor\begin{matrix}
A_1,A_2,\dots,A_{r+1}\\C_1,C_2,\dots,C_r
\end{matrix};Q,Z\right\rceil:=
\sum_{k\ge 0}
\left\lfloor\begin{matrix}
A_1,A_2,\dots,A_{r+1}\\C_1,C_2,\dots,C_r,Q
\end{matrix};Q,Z\right\rceil_k.
\end{equation}
We also refer to the respective series as
{\em (noncommutative) $Q$-hypergeometric series}.
In each case, the series terminates if one of the upper
parameters $A_i$ is of the form $Q^{-n}$. If the series
does not terminate, then (implicitly assuming that
$R$ is Banach algebra with norm $\|\cdot\|$) it converges
if $\|Z\|<1$.

Note that the factors in the generalized noncommutative
($Q$-)shifted factorials are strongly interlaced, e.g.,
\begin{equation*}
\left\lfloor\begin{matrix}A,B\\C,D\end{matrix};Z\right\rceil_2=
C^{-1}AD^{-1}BZ(C+I)^{-1}(A+I)(D+I)^{-1}(B+I)Z.
\end{equation*}
It is this interlacing which is mainly responsible that the
noncommutative hypergeometric series considered in this paper
can be summed in closed form. This is maybe best understood by
regarding the general procedure for proving (all) the terminating
identities in this paper, namely induction:
A particular factor of the summand is usually
rewritten such that the original sum is split into two sums.
After shifting the index of summation in one of the sums some
factors can be pulled out and a similar sum remains to which the
inductive hypothesis applies. (See Sections~\ref{sec2f1}--\ref{sec7f6}
for several demonstrations of this procedure.) The ($Q$-)shifted
factorials of types I and II have been defined exactly in a way that
induction can be successfully applied for proving the respective
summations.

\begin{Remark}
Tirao's~\cite{T} matrix extension of the Gau{\ss} hypergeometric
function corresponds to the special case of our noncommutative
${}_2F_1$ series of type I when the parameters are $n\times n$
matrices over the complex numbers $\mathbb C$ and, in addition,
the argument $Z$ is a diagonal matrix $zI$ with $z\in\mathbb C$.
Restated in terms of the notation introduced in this section,
Tirao essentially shows (among other results) that
\cite[Th.~2]{T} {\it if $A,B,C,F_0\in R$ and $C+jI$ is invertible
for all nonnegative integers $j$, then
\begin{equation}\label{2f1typeI}
F(z)={}_2F_1\!\left\lceil\begin{matrix}
A,B\\C\end{matrix};zI\right\rfloor\,F_0
\end{equation}
is analytic on $|z|<1$ with values in $R$, and $F(z)$ is a
solution of the hypergeometric equation
\begin{equation}
z(1-z)F''+[C-z(I+A+B)]F'-ABF=O
\end{equation}
such that $F(0)=F_0$, and conversely any solution of $F$
analytic at $z=0$ is of this form.} 
He further shows (see \cite[Cor.~3]{T}) that the
matrix valued Jacobi polynomials introduced by Gr\"unbaum~\cite{G}
can be expressed in terms of hypergeometric functions of the above type
\eqref{2f1typeI}, thereby giving an explicit example within the theory
of matrix valued orthogonal polynomials initiated by Krein~\cite{Kr}.

In view of the above, it appears indeed very appropriate to study the
${}_2F_1$ function of type II in \eqref{defnchypII} and even the
the ${}_2\phi_1$ functions of types I and II, in \eqref{defnchypIQ}
and \eqref{defnchypIIQ}, respectively, in terms of the second-order
($Q$-)differential equations they (most likely) satisfy. We defer the
investigation of this interesting question to elsewhere
as it goes beyond the scope of this paper, the (elementary)
derivation of explicit summation formulae.
\end{Remark}

\section{Elementary identities for noncommutative ($Q$-)shifted factorials}
\label{secelid}

Here we provide a couple of lemmas which will be utilized for proving
the summations in Sections~\ref{sec2f1}--\ref{sec7f6}. These lemmas
concern addition formulae for noncommutative ($Q$-)shifted factorials
of type I and type II. Throughout we assume $n$ to be a nonnegative
integer. All formulae are proved in the same manner,
by induction on $n$.

\begin{Lemma}\label{lem21}
Let $A$ and $C$ be noncommutative parameters of some unit ring,
and let $n$ be a nonnegative integer.
Then we have the following addition formula
for shifted factorials of type I.
\begin{equation}\label{lem21Igl}
\left\lceil\begin{matrix}C-A\\C\end{matrix}\right\rfloor_n
-\left\lceil\begin{matrix}C-A\\C+I\end{matrix}\right\rfloor_n
C^{-1}A=
\left\lceil\begin{matrix}C-A\\C\end{matrix}\right\rfloor_{n+1}.
\end{equation}
Further, we have the following addition formula
for shifted factorials of type II.
\begin{align}\label{lem21IIgl}\notag
\left\lfloor\begin{matrix}C-A+I\\C\end{matrix}\right\rceil_n
(C-A+nI)^{-1}-
C^{-1}A\left\lfloor\begin{matrix}C-A+I\\C+I\end{matrix}\right\rceil_n
(C-A+nI)^{-1}&\\=
\left\lfloor\begin{matrix}C-A+I\\C\end{matrix}\right\rceil_{n+1}
(C-A+(n+1)I)^{-1}&.
\end{align}
\end{Lemma}

\begin{proof}
We start with \eqref{lem21Igl}.
For $n=0$ \eqref{lem21Igl} is easily verified.
Assume the identity is true for all nonnegative integers
less than a fixed positive integer $n$.
Then we rewrite the left-hand side of \eqref{lem21Igl} as
\begin{multline*}
\left\lceil\begin{matrix}C-A\\C\end{matrix}\right\rfloor_n
-\left\lceil\begin{matrix}C-A\\C+I\end{matrix}\right\rfloor_n
C^{-1}A\\=
\left\lceil\begin{matrix}C-A+I\\C+I\end{matrix}\right\rfloor_{n-1}
C^{-1}(C-A)
-\left\lceil\begin{matrix}C-A\\C+I\end{matrix}\right\rfloor_n
C^{-1}A.
\end{multline*}
We have simply pulled out the last two factors from the first term.
Now we can apply the inductive hypothesis (with $C$ replaced by $C+I$)
to the first term to transform the last expression into
\begin{align*}
\left(\left\lceil\begin{matrix}C-A+I\\C+I\end{matrix}\right\rfloor_n
+\left\lceil\begin{matrix}C-A+I\\C+2I\end{matrix}\right\rfloor_{n-1}
(C+I)^{-1}A\right)C^{-1}(C-A)
-\left\lceil\begin{matrix}C-A\\C+I\end{matrix}\right\rfloor_n
C^{-1}A&\\=
\left\lceil\begin{matrix}C-A\\C\end{matrix}\right\rfloor_{n+1}
+\left\lceil\begin{matrix}C-A+I\\C+2I\end{matrix}\right\rfloor_{n-1}
(C+I)^{-1}AC^{-1}(C-A)
-\left\lceil\begin{matrix}C-A\\C+I\end{matrix}\right\rfloor_n
C^{-1}A&.
\end{align*}
What remains to be shown is that in this sum of three terms
the last two cancel each other. This is equivalent to
\begin{equation}\label{easyidI}
AC^{-1}(C-A)=(C-A)C^{-1}A,
\end{equation}
which is immediately verified since both sides equal $A-AC^{-1}A$.
Hence, we have established \eqref{lem21Igl}.

Next, we prove \eqref{lem21IIgl}. The $n=0$ case is trivial.
Next, assume that the formula is true for all nonnegative integers
less than a fixed positive integer $n$.
Then we rewrite the left-hand side of \eqref{lem21IIgl} as
\begin{align*}
\left\lfloor\begin{matrix}C-A+I\\C\end{matrix}\right\rceil_n
(C-A+nI)^{-1}-
C^{-1}A\left\lfloor\begin{matrix}C-A+I\\C+I\end{matrix}\right\rceil_n
(C-A+nI)^{-1}&\\=
C^{-1}(C-A+I)\left\lfloor
\begin{matrix}C-A+2I\\C+I\end{matrix}\right\rceil_{n-1}
(C-A+nI)^{-1}&\\
-C^{-1}A\left\lfloor\begin{matrix}C-A+I\\C+I\end{matrix}\right\rceil_n
(C-A+nI)^{-1}&.
\end{align*}
We have simply pulled out the first two factors from the first term.
Now we can apply the inductive hypothesis (with $C$ replaced by $C+I$)
to the first term to transform the last expression into
\begin{multline*}
C^{-1}(C-A+I)
\bigg(\left\lfloor\begin{matrix}C-A+2I\\C+I\end{matrix}\right\rceil_n
(C-A+(n+1)I)^{-1}\\
+(C+I)^{-1}A\left\lfloor
\begin{matrix}C-A+2I\\C+2I\end{matrix}\right\rceil_{n-1}
(C-A+nI)^{-1}\bigg)\\
-C^{-1}A\left\lfloor\begin{matrix}C-A+I\\C+I\end{matrix}\right\rceil_n
(C-A+nI)^{-1}.
\end{multline*}
What remains to be shown is that in the resulting sum of three terms
the last two cancel each other, which is equivalent to
\begin{equation*}
(C-A+I)(C+I)^{-1}A=
A(C+I)^{-1}(C-A+I).
\end{equation*}
But this is simply the $C\mapsto C+I$ case of \eqref{easyidI}.
Thus, we have established \eqref{lem21IIgl}.
\end{proof}

\begin{Lemma}\label{lem21Q}
Let $A$ and $C$ be noncommutative parameters of some unit ring,
and suppose that $Q$ commutes both with $A$ and $C$. Further,
let $n$ be a nonnegative integer.
Then we have the following addition formulae
for $Q$-shifted factorials of type I.
\begin{equation}\label{lem21QIgl}
\left\lceil\begin{matrix}CA^{-1}\\C\end{matrix};Q,A\right\rfloor_n
-\left\lceil\begin{matrix}CA^{-1}\\CQ\end{matrix};Q,A\right\rfloor_n
(I-C)^{-1}(I-A)=
\left\lceil\begin{matrix}CA^{-1}\\C\end{matrix};Q,A\right\rfloor_{n+1},
\end{equation}
and
\begin{equation}\label{lem21QIrgl}
\left\lceil\begin{matrix}A^{-1}C\\C\end{matrix};Q,I\right\rfloor_n
-\left\lceil\begin{matrix}A^{-1}C\\CQ\end{matrix};Q,I\right\rfloor_n
(I-C)^{-1}(I-A)A^{-1}CQ^n=
\left\lceil\begin{matrix}A^{-1}C\\C\end{matrix};Q,I\right\rfloor_{n+1}.
\end{equation}
Further, we have the following addition formulae
for $Q$-shifted factorials of type II.
\begin{multline}\label{lem21QIIgl}
\left\lfloor\begin{matrix}CA^{-1}Q\\C\end{matrix};Q,A\right\rceil_n
(A-CQ^n)^{-1}-(I-C)^{-1}(I-A)
\left\lfloor\begin{matrix}CA^{-1}Q\\CQ\end{matrix};Q,A\right\rceil_n
(A-CQ^n)^{-1}\\=
\left\lfloor\begin{matrix}CA^{-1}Q\\C\end{matrix};Q,A\right\rceil_{n+1}
(A-CQ^{n+1})^{-1},
\end{multline}
and
\begin{multline}\label{lem21QIIrgl}
\left\lfloor\begin{matrix}A^{-1}CQ\\C\end{matrix};Q,I\right\rceil_n
(I-A^{-1}CQ^n)^{-1}\\
-(I-C)^{-1}(I-A)A^{-1}CQ^n
\left\lfloor\begin{matrix}A^{-1}CQ\\CQ\end{matrix};Q,I\right\rceil_n
(I-A^{-1}CQ^n)^{-1}\\=
\left\lfloor\begin{matrix}A^{-1}CQ\\C\end{matrix};Q,I\right\rceil_{n+1}
(I-A^{-1}CQ^{n+1})^{-1}.
\end{multline}
\end{Lemma}

\begin{proof}
We restrict ourselves to proving \eqref{lem21QIIgl}, the proofs
of \eqref{lem21QIgl}, \eqref{lem21QIrgl} and \eqref{lem21QIIrgl}
being similar. For $n=0$ \eqref{lem21QIIgl} is easily verified.
Assume the identity is true for all nonnegative integers
less than a fixed positive integer $n$.
Then we rewrite the left-hand side of \eqref{lem21QIIgl} as
\begin{align*}
\!\left\lfloor\begin{matrix}CA^{-1}Q\\C\end{matrix};Q,A\right\rceil_n
(A-CQ^n)^{-1}-(I-C)^{-1}(I-A)
\left\lfloor\begin{matrix}CA^{-1}Q\\CQ\end{matrix};Q,A\right\rceil_n
(A-CQ^n)^{-1}&\\=
(I-C)^{-1}(A-CQ)\left\lfloor
\begin{matrix}CA^{-1}Q^2\\CQ\end{matrix};Q,A\right\rceil_{n-1}
(A-CQ^n)^{-1}&\\
-(I-C)^{-1}(I-A)\left\lfloor\begin{matrix}CA^{-1}Q\\
CQ\end{matrix};Q,A\right\rceil_n(A-CQ^n)^{-1}.&
\end{align*}
We have simply pulled out the first two factors from the first term.
Now we can apply the inductive hypothesis (with $C$ replaced by $CQ$)
to the first term to transform the last expression into
\begin{multline*}
(I-C)^{-1}(A-CQ)
\left(\left\lfloor\begin{matrix}CA^{-1}Q^2\\CQ\end{matrix};Q,A\right\rceil_n
(A-CQ^{n+1})^{-1}\right.\\
\left.+(I-CQ)^{-1}(I-A)\left\lfloor
\begin{matrix}CA^{-1}Q^2\\CQ^2\end{matrix};Q,A\right\rceil_{n-1}
(A-CQ^n)^{-1}\right)\\
-(I-C)^{-1}(I-A)\left\lfloor\begin{matrix}CA^{-1}Q\\
CQ\end{matrix};Q,A\right\rceil_n
(A-CQ^n)^{-1}.
\end{multline*}
What remains to be shown is that in the resulting sum of three terms
the last two cancel each other, which is equivalent to
\begin{equation}\label{easyidQ}
(A-CQ)(I-CQ)^{-1}(I-A)=
(I-A)(I-CQ)^{-1}(A-CQ).
\end{equation}
However, splitting the factor $(A-CQ)$ on each side of \eqref{easyidQ}
into two terms as $(I-CQ)-(I-A)$, both sides can be reduced to the same
expression, namely
\begin{equation*}
(I-A)-(I-A)(I-CQ)^{-1}(I-A),
\end{equation*}
which immediately establishes \eqref{lem21QIIgl}
\end{proof}

\begin{Lemma}\label{lem32}
Let $A$, $B$ and $C$ be noncommutative parameters of some unit ring,
and suppose that the sum $A+B-C$ commutes each with $A$, $B$ and $C$.
Further, let $n$ be a nonnegative integer. Then
we have the following addition formula
for shifted factorials of type I.
\begin{multline}\label{lem32Igl}
\left\lceil\begin{matrix}C-B,C-A\\
C,C-A-B\end{matrix}\right\rceil_n-
\left\lfloor\begin{matrix}C-B,C-A\\
C+I,C-A-B-I\end{matrix}\right\rceil_n\\\times
C^{-1}A(A+B-C+(1-n)I)^{-1}B(A+B-C-nI)^{-1}(A+B-C+I)\\=
\left\lfloor\begin{matrix}C-B,C-A\\
C,C-A-B\end{matrix}\right\rceil_{n+1}.
\end{multline}
Further, we have the following addition formula
for shifted factorials of type II.
\begin{multline}\label{lem32IIgl}
\left\lfloor\begin{matrix}C-B+I,C-A+I\\
C,C-A-B\end{matrix}\right\rceil_n
(C-A+nI)^{-1}(C-B+nI)^{-1}\\
-C^{-1}A(A+B-C+(1-n)I)^{-1}B(A+B-C-nI)^{-1}(A+B-C+I)\\\times
\left\lfloor\begin{matrix}C-B+I,C-A+I\\
C+I,C-A-B-I\end{matrix}\right\rceil_n
(C-A+nI)^{-1}(C-B+nI)^{-1}\\=
\left\lfloor\begin{matrix}C-B+I,C-A+I\\
C,C-A-B\end{matrix}\right\rceil_{n+1}
(C-A+(n+1)I)^{-1}(C-B+(n+1)I)^{-1}.
\end{multline}
\end{Lemma}

\begin{proof}
We prove \eqref{lem32IIgl} by induction on $n$ and
leave the proof of \eqref{lem32Igl} (which is similar) to the reader.
For $n=0$ \eqref{lem32IIgl} is easily verified.
Assume the identity is true for all nonnegative integers
less than a fixed positive integer $n$.
Then we rewrite the left-hand side of \eqref{lem32IIgl} as
\begin{multline*}
\left\lfloor\begin{matrix}C-B+I,C-A+I\\
C,C-A-B\end{matrix}\right\rceil_n
(C-A+nI)^{-1}(C-B+nI)^{-1}\\
-C^{-1}A(A+B-C+(1-n)I)^{-1}B(A+B-C-nI)^{-1}(A+B-C+I)\\\times
\left\lfloor\begin{matrix}C-B+I,C-A+I\\
C+I,C-A-B-I\end{matrix}\right\rceil_n
(C-A+nI)^{-1}(C-B+nI)^{-1}\\=
C^{-1}(C-B+I)(C-A-B)^{-1}(C-A+I)\\\times
\left\lfloor\begin{matrix}C-B+2I,C-A+2I\\
C+I,C-A-B+I\end{matrix}\right\rceil_{n-1}
(C-A+nI)^{-1}(C-B+nI)^{-1}\\
-C^{-1}A(A+B-C+(1-n)I)^{-1}B(A+B-C-nI)^{-1}(A+B-C+I)\\\times
\left\lfloor\begin{matrix}C-B+I,C-A+I\\
C+I,C-A-B-I\end{matrix}\right\rceil_n
(C-A+nI)^{-1}(C-B+nI)^{-1}.
\end{multline*}
We have simply pulled out the first four factors from the first term.
Now we can apply the inductive hypothesis (with $C$ replaced by $C+I$)
to the first term to transform the last expression into
\begin{multline*}
C^{-1}(C-B+I)(C-A-B)^{-1}(C-A+I)\\\times
\left(\left\lfloor\begin{matrix}C-B+2I,C-A+2I\\
C+I,C-A-B+I\end{matrix}\right\rceil_n
(C-A+(n+1)I)^{-1}(C-B+(n+1)I)^{-1}\right.\\
+(C+I)^{-1}A(A+B-C+(1-n)I)^{-1}B(A+B-C-nI)^{-1}(A+B-C)\\\times
\left.\left\lfloor\begin{matrix}C-B+2I,C-A+2I\\
C+2I,C-A-B\end{matrix}\right\rceil_{n-1}
(C-A+nI)^{-1}(C-B+nI)^{-1}\right)\\
-C^{-1}A(A+B-C+(1-n)I)^{-1}B(A+B-C-nI)^{-1}(A+B-C+I)\\\times
\left\lfloor\begin{matrix}C-B+I,C-A+I\\
C+I,C-A-B-I\end{matrix}\right\rceil_n
(C-A+nI)^{-1}(C-B+nI)^{-1}.
\end{multline*}
What remains to be shown is that in the resulting sum of three terms
the last two cancel each other, which is equivalent to
\begin{equation}\label{easy32id}
(C-B+I)(C-A+I)(C+I)^{-1}AB=
AB(C+I)^{-1}(C-B+I)(C-A+I).
\end{equation}
However, splitting the factor $(C-B+I)(C-A+I)$ on each side of
\eqref{easy32id} into two terms as $AB+(C-A-B+I)(C+I)$
(which can be done since $A(C-B)=(C-B)A$),
both sides can be reduced to the same
expression, namely
\begin{equation*}
AB(C+I)^{-1}AB+(C-A-B+I)AB,
\end{equation*}
which immediately establishes \eqref{lem32IIgl}.
\end{proof}

\begin{Lemma}\label{lem32Q}
Let $A$, $B$ and $C$ be noncommutative parameters of some unit ring,
and suppose that $Q$ commutes each with $A$, $B$ and $C$. Further,
assume that the product $BC^{-1}A$ commutes each with $A$, $B$ and $C$.
Moreover, let $n$ be a nonnegative integer. Then
we have the following addition formula
for $Q$-shifted factorials of type I.
\begin{multline}\label{lem32QIgl}
\left\lfloor\begin{matrix}CB^{-1},A^{-1}C\\
C,A^{-1}CB^{-1}\end{matrix};Q,I\right\rceil_n
-\left\lfloor\begin{matrix}CB^{-1},A^{-1}C\\
CQ,A^{-1}CB^{-1}Q^{-1}\end{matrix};Q,I\right\rceil_n\\\times
(I-C)^{-1}(I-A)(I-BC^{-1}AQ^{1-n})^{-1}(I-B)
(I-BC^{-1}AQ^{-n})^{-1}(I-BC^{-1}AQ)\\=
\left\lfloor\begin{matrix}CB^{-1},A^{-1}C\\
C,A^{-1}CB^{-1}\end{matrix};Q,I\right\rceil_{n+1}.
\end{multline}
Further, we have the following addition formula
for $Q$-shifted factorials of type II.
\begin{multline}\label{lem32QIIgl}
\left\lfloor\begin{matrix}CB^{-1}Q,A^{-1}CQ\\
C,A^{-1}CB^{-1}\end{matrix};Q,I\right\rceil_n
(I-A^{-1}CQ^n)^{-1}(I-CB^{-1}Q^n)^{-1}\\
-(I-C)^{-1}(I-A)(I-BC^{-1}AQ^{1-n})^{-1}(I-B)
(I-BC^{-1}AQ^{-n})^{-1}(I-BC^{-1}AQ)\\\times
\left\lfloor\begin{matrix}CB^{-1}Q,A^{-1}CQ\\
CQ,A^{-1}CB^{-1}Q^{-1}\end{matrix};Q,I\right\rceil_n
(I-A^{-1}CQ^n)^{-1}(I-CB^{-1}Q^n)^{-1}\\=
\left\lfloor\begin{matrix}CB^{-1}Q,A^{-1}CQ\\
C,A^{-1}CB^{-1}\end{matrix};Q,I\right\rceil_{n+1}
(I-A^{-1}CQ^{n+1})^{-1}(I-CB^{-1}Q^{n+1})^{-1}.
\end{multline}
\end{Lemma}
\begin{proof}
We prove \eqref{lem32QIIgl} by induction on $n$ and
leave the proof of \eqref{lem32QIgl} (which is similar) to the reader.
For $n=0$ \eqref{lem32QIIgl} is easily verified.
Assume the identity is true for all nonnegative integers
less than a fixed positive integer $n$.
Then we rewrite the left-hand side of \eqref{lem32QIIgl} as
\begin{multline*}
\left\lfloor\begin{matrix}CB^{-1}Q,A^{-1}CQ\\
C,A^{-1}CB^{-1}\end{matrix};Q,I\right\rceil_n
(I-A^{-1}CQ^n)^{-1}(I-CB^{-1}Q^n)^{-1}\\
-(I-C)^{-1}(I-A)(I-BC^{-1}AQ^{1-n})^{-1}(I-B)
(I-BC^{-1}AQ^{-n})^{-1}(I-BC^{-1}AQ)\\\times
\left\lfloor\begin{matrix}CB^{-1}Q,A^{-1}CQ\\
CQ,A^{-1}CB^{-1}Q^{-1}\end{matrix};Q,I\right\rceil_n
(I-A^{-1}CQ^n)^{-1}(I-CB^{-1}Q^n)^{-1}\\=
(I-C)^{-1}(I-CB^{-1}Q)(I-A^{-1}CB^{-1})^{-1}(I-A^{-1}CQ)\\\times
\left\lfloor\begin{matrix}CB^{-1}Q^2,A^{-1}CQ^2\\
CQ,A^{-1}CB^{-1}Q\end{matrix};Q,I\right\rceil_{n-1}
(I-A^{-1}CQ^n)^{-1}(I-CB^{-1}Q^n)^{-1}\\
-(I-C)^{-1}(I-A)(I-BC^{-1}AQ^{1-n})^{-1}(I-B)
(I-BC^{-1}AQ^{-n})^{-1}(I-BC^{-1}AQ)\\\times
\left\lfloor\begin{matrix}CB^{-1}Q,A^{-1}CQ\\
CQ,A^{-1}CB^{-1}Q^{-1}\end{matrix};Q,I\right\rceil_n
(I-A^{-1}CQ^n)^{-1}(I-CB^{-1}Q^n)^{-1}.
\end{multline*}
We have simply pulled out the first four factors from the first term.
Now we can apply the inductive hypothesis (with $C$ replaced by $CQ$)
to the first term to transform the last expression into
\begin{multline*}
(I-C)^{-1}(I-CB^{-1}Q)(I-A^{-1}CB^{-1})^{-1}(I-A^{-1}CQ)\\\times
\left(\left\lfloor\begin{matrix}CB^{-1}Q^2,A^{-1}CQ^2\\
CQ,A^{-1}CB^{-1}Q^{-1}\end{matrix};Q,I\right\rceil_n
(I-A^{-1}CQ^{n+1})^{-1}(I-CB^{-1}Q^{n+1})^{-1}\right.\\
+(I-CQ)^{-1}(I-A)(I-BC^{-1}AQ^{1-n})^{-1}(I-B)
(I-BC^{-1}AQ^{-n})^{-1}(I-BC^{-1}A)\\\times
\left.\left\lfloor\begin{matrix}CB^{-1}Q^2,A^{-1}CQ^2\\
CQ^2,A^{-1}CB^{-1}\end{matrix};Q,I\right\rceil_{n-1}
(I-A^{-1}CQ^n)^{-1}(I-CB^{-1}Q^n)^{-1}\right)\\
-(I-C)^{-1}(I-A)(I-BC^{-1}AQ^{1-n})^{-1}(I-B)
(I-BC^{-1}AQ^{-n})^{-1}(I-BC^{-1}AQ)\\\times
\left\lfloor\begin{matrix}CB^{-1}Q,A^{-1}CQ\\
CQ,A^{-1}CB^{-1}Q^{-1}\end{matrix};Q,I\right\rceil_n
(I-A^{-1}CQ^n)^{-1}(I-CB^{-1}Q^n)^{-1}.
\end{multline*}
What remains to be shown is that in the resulting sum of three terms
the last two cancel each other, which is equivalent to
\begin{multline}\label{easy32Qid}
(I-CB^{-1}Q)(I-A^{-1}CQ)(I-CQ)^{-1}(I-A)(I-B)\\=
(I-A)(I-B)(I-CQ)^{-1}(I-CB^{-1}Q)(I-A^{-1}CQ).
\end{multline}
However, splitting the factor $(I-CB^{-1}Q)(I-A^{-1}CQ)$ on each side of
\eqref{easy32Qid} into two terms as
$(I-A)(I-B)A^{-1}CB^{-1}Q+(I-A^{-1}CB^{-1}Q)(I-CQ)$
(which can be done since $CB^{-1}A^{-1}=A^{-1}CB^{-1}$, etc.),
both sides can be reduced to the same
expression, namely
\begin{multline*}
(I-A)(I-B)(I-CQ)^{-1}(I-A)(I-B)A^{-1}CB^{-1}Q\\
+(I-A^{-1}CB^{-1}Q)(I-A)(I-B),
\end{multline*}
which immediately establishes \eqref{lem32QIIgl}.
\end{proof}

\begin{Lemma}\label{lem76}
Let $A$, $B$, $C$ and $D$ be noncommutative parameters of some unit ring.
Assume that the commutation relations \eqref{comrel} hold.
Further, let $n$ be a nonnegative integer. Then
we have the following addition formula
for shifted factorials of type I.
\begin{multline}\label{lem76gl}
\left\lceil\begin{matrix}A-C-D+I,A-B-D+I,A+I,A-B-C+I\\
A-C+I,A-B+I,A-D+I,A-B-C-D+I\end{matrix}\right\rfloor_n\\-
\left\lceil\begin{matrix}A-C-D+I,A-B-D+I,A+3I,A-B-C+I\\
A-C+2I,A-B+2I,A-D+2I,A-B-C-D\end{matrix}\right\rfloor_n\\\times
(B+C+D-A-(n+1)I)^{-1}(B+C+D-A)(A+(n+1)I)^{-1}(A+I)\\\times
(A+(n+2)I)^{-1}(A+2I)(B+C+D-A-nI)^{-1}(2A-B-C-D+(2+2n)I)\\\times
(A-C+I)^{-1}B(A-B+I)^{-1}C(A-D+I)^{-1}D\\=
\left\lceil\begin{matrix}A-C-D+I,A-B-D+I,A+I,A-B-C+I\\
A-C+I,A-B+I,A-D+I,A-B-C-D+I\end{matrix}\right\rfloor_{n+1}.
\end{multline}
\end{Lemma}
\begin{Remark}
An equivalent, almost identical formula holds for shifted factorials of
type II, see also Remark~\ref{remIaII}.
\end{Remark}
\begin{proof}[Proof of Lemma~\ref{lem76}]
We prove \eqref{lem76gl} by induction on $n$.
For $n=0$ \eqref{lem76gl} is easily verified.
Assume the identity is true for all nonnegative integers
less than a fixed positive integer $n$.
Then we rewrite the left-hand side of \eqref{lem76gl} as
\begin{multline*}
\left\lceil\begin{matrix}A-C-D+I,A-B-D+I,A+I,A-B-C+I\\
A-C+I,A-B+I,A-D+I,A-B-C-D+I\end{matrix}\right\rfloor_n\\-
\left\lceil\begin{matrix}A-C-D+I,A-B-D+I,A+3I,A-B-C+I\\
A-C+2I,A-B+2I,A-D+2I,A-B-C-D\end{matrix}\right\rfloor_n\\\times
(B+C+D-A-(n+1)I)^{-1}(B+C+D-A)(A+(n+1)I)^{-1}(A+I)\\\times
(A+(n+2)I)^{-1}(A+2I)(B+C+D-A-nI)^{-1}(2A-B-C-D+(2+2n)I)\\\times
(A-C+I)^{-1}B(A-B+I)^{-1}C(A-D+I)^{-1}D\\=
\left\lceil\begin{matrix}A-C-D+2I,A-B-D+2I,A+2I,A-B-C+2I\\
A-C+2I,A-B+2I,A-D+2I,A-B-C-D+2I\end{matrix}\right\rfloor_{n-1}\\\times
(A-C+I)^{-1}(A-C-D+I)(A-B+I)^{-1}(A-B-D+I)\\\times
(A-D+I)^{-1}(A+I)(A-B-C-D+I)^{-1}(A-B-C+I)\\-
\left\lceil\begin{matrix}A-C-D+I,A-B-D+I,A+3I,A-B-C+I\\
A-C+2I,A-B+2I,A-D+2I,A-B-C-D\end{matrix}\right\rfloor_n\\\times
(B+C+D-A-(n+1)I)^{-1}(B+C+D-A)(A+(n+1)I)^{-1}(A+I)\\\times
(A+(n+2)I)^{-1}(A+2I)(B+C+D-A-nI)^{-1}(2A-B-C-D+(2+2n)I)\\\times
(A-C+I)^{-1}B(A-B+I)^{-1}C(A-D+I)^{-1}D.
\end{multline*}
We have simply pulled out the first eight factors from the first term.
Now we can apply the inductive hypothesis (with $A$ replaced by $A+I$)
to the first term to transform the last expression into
\begin{multline*}
\bigg(\left\lceil\begin{matrix}A-C-D+2I,A-B-D+2I,A+2I,A-B-C+2I\\
A-C+2I,A-B+2I,A-D+2I,A-B-C-D+2I\end{matrix}\right\rfloor_n\\+
\left\lceil\begin{matrix}A-C-D+2I,A-B-D+2I,A+4I,A-B-C+2I\\
A-C+3I,A-B+3I,A-D+3I,A-B-C-D+I\end{matrix}\right\rfloor_{n-1}\\\times
(B+C+D-A-(n+1)I)^{-1}(B+C+D-A-I)(A+(n+1)I)^{-1}(A+2I)\\\times
(A+(n+2)I)^{-1}(A+3I)(B+C+D-A-nI)^{-1}(2A-B-C-D+(2+2n)I)\\\times
(A-C+2I)^{-1}B(A-B+2I)^{-1}C(A-D+2I)^{-1}D\bigg)\\\times
(A-C+I)^{-1}(A-C-D+I)(A-B+I)^{-1}(A-B-D+I)\\\times
(A-D+I)^{-1}(A+I)(A-B-C-D+I)^{-1}(A-B-C+I)\\
-\left\lceil\begin{matrix}A-C-D+I,A-B-D+I,A+3I,A-B-C+I\\
A-C+2I,A-B+2I,A-D+2I,A-B-C-D\end{matrix}\right\rfloor_n\\\times
(B+C+D-A-(n+1)I)^{-1}(B+C+D-A)(A+(n+1)I)^{-1}(A+I)\\\times
(A+(n+2)I)^{-1}(A+2I)(B+C+D-A-nI)^{-1}(2A-B-C-D+(2+2n)I)\\\times
(A-C+I)^{-1}B(A-B+I)^{-1}C(A-D+I)^{-1}D.
\end{multline*}
What remains to be shown is that in the resulting sum of three terms
the last two cancel each other, which is equivalent to
\begin{multline*}
BCD(A-D+2I)^{-1}(A-C+I)^{-1}(A-B+I)^{-1}\\\times
(A-C-D+I)(A-B-D+I)(A-B-C+I)\\
=(A-C-D+I)(A-B-D+I)(A-B-C+I)\\\times
(A-D+2I)^{-1}(A-C+I)^{-1}(A-B+I)^{-1}BCD.
\end{multline*}
This follows from the fact that the three products $BCD$,
$(A-B+I)(A-C+I)(A-D+2I)$ and $(A-C-D+I)(A-B-D+I)(A-B-C+I)$ mutually
commute, which can be readily verified using \eqref{comrel} and
\eqref{comrelbcdd}.
\end{proof}

\section{Chu--Vandermonde summations}\label{sec2f1}

The classical Chu--Vandermonde summation formula
(cf.\ \cite[Appendix~(III.4)]{S})
sums a terminating ${}_2F_1$ series with unit argument:
\begin{equation}\label{2f1}
{}_2F_1\!\left[\begin{matrix}a,-n\\c\end{matrix};1\right]=
\frac{(c-a)_n}{(c)_n}.
\end{equation}

We provide two noncommutative extensions of \eqref{2f1}:
\begin{Theorem}\label{nc2f1}
Let $A$, $C$ be noncommutative parameters of some unit ring,
and let $n$ be a nonnegative integer. Then we have the following
summation for a noncommutative hypergeometric series of type I.
\begin{equation}\label{nc2f1Igl}
{}_2F_1\!\left\lceil\begin{matrix}A,-nI\\C\end{matrix};I\right\rfloor=
\left\lceil\begin{matrix}C-A\\C\end{matrix}\right\rfloor_n.
\end{equation}
Further, we have the following
summation for a noncommutative hypergeometric series of type II.
\begin{equation}\label{nc2f1IIgl}
{}_2F_1\!\left\lfloor\begin{matrix}A,-nI\\C\end{matrix};I\right\rceil=
\left\lfloor\begin{matrix}C-A+I\\C\end{matrix}\right\rceil_n
(C-A+nI)^{-1}(C-A).
\end{equation}
\end{Theorem}
Note that the right-hand side of \eqref{nc2f1IIgl} may also be written as
\begin{equation*}
\left\lfloor\begin{matrix}C-A+I\\C\end{matrix}\right\rceil_{n-1}
(C+(n-1)I)^{-1}(C-A),
\end{equation*}
however, in the simple case $n=0$ it is a little bit easier to
see that the sum reduces correctly to $I$ when one uses the expression
on the right-hand side of \eqref{nc2f1IIgl}.

\begin{proof}[Proof of Theorem~\ref{nc2f1}]
Both identities \eqref{nc2f1Igl} and \eqref{nc2f1IIgl}
are readily proved by induction on $n$. We start with
the proof of the first summation \eqref{nc2f1Igl}.

For $n=0$ the formula is trivial.
We assume the summation is true up to a fixed $n$.
To prove \eqref{nc2f1Igl} for $n+1$, use the elementary identity
\begin{equation*}
\left\lceil\begin{matrix}-(n+1)I\\I\end{matrix}\right\rfloor_k=
\left\lceil\begin{matrix}-nI\\I\end{matrix}\right\rfloor_k
\left[I-((-n-1+k)I)^{-1}kI\right]
\end{equation*}
to obtain
\begin{align*}
{}_2F_1\!\left\lceil\begin{matrix}A,-(n+1)I\\C\end{matrix};I\right\rfloor=
{}&\sum_{k=0}^n
\left\lceil\begin{matrix}A,-nI\\C,I\end{matrix}\right\rfloor_k
\left[I-((-n-1+k)I)^{-1}kI\right]\\=
{}&{}_2F_1\!\left\lceil\begin{matrix}A,-nI\\C\end{matrix};I\right\rfloor-
{}_2F_1\!\left\lceil\begin{matrix}A+I,-nI\\C+I\end{matrix};I\right\rfloor
\;C^{-1}A\\=
{}&\left\lceil\begin{matrix}C-A\\C\end{matrix}\right\rfloor_n-
\left\lfloor\begin{matrix}C-A\\C+I\end{matrix}\right\rceil_n
C^{-1}A=
\left\lfloor\begin{matrix}C-A\\C\end{matrix}\right\rceil_{n+1},
\end{align*}
the penultimate equation due to the inductive hypothesis,
the last equation due to Lemma~\ref{lem21}, Equation~\eqref{lem21Igl}.

We turn now to the proof of the second identity \eqref{nc2f1IIgl}.
For $n=0$ the formula is trivial.
We assume the summation is true up to a fixed $n$.
To prove \eqref{nc2f1IIgl} for $n+1$, use the elementary identity
\begin{equation*}
\left\lfloor\begin{matrix}-(n+1)I\\I\end{matrix}\right\rceil_k=
\left\lfloor\begin{matrix}-nI\\I\end{matrix}\right\rceil_k
\left[I-((-n-1+k)I)^{-1}kI\right]
\end{equation*}
to obtain
\begin{align*}
{}_2F_1\!\left\lfloor\begin{matrix}A,-(n+1)I\\C\end{matrix};I\right\rceil=
{}&\sum_{k=0}^n
\left\lfloor\begin{matrix}A,-nI\\C,I\end{matrix}\right\rceil_k
\left[I-((-n-1+k)I)^{-1}kI\right]\\=
{}&{}_2F_1\!\left\lfloor\begin{matrix}A,-nI\\C\end{matrix};I\right\rceil-
C^{-1}A\;{}_2F_1\!\left\lfloor
\begin{matrix}A+I,-nI\\C+I\end{matrix};I\right\rceil\\=
{}&\left\lfloor\begin{matrix}C-A+I\\C\end{matrix}\right\rceil_n
(C-A+nI)^{-1}(C-A)\\&-
C^{-1}A\left\lfloor\begin{matrix}C-A+I\\C+I\end{matrix}\right\rceil_n
(C-A+nI)^{-1}(C-A)\\=
{}&\left\lfloor\begin{matrix}C-A+I\\C\end{matrix}\right\rceil_{n+1}
(C-A+(n+1)I)^{-1}(C-A),
\end{align*}
the penultimate equation due to the inductive hypothesis,
the last equation due to Lemma~\ref{lem21}, Equation~\eqref{lem21IIgl}.
\end{proof}

The following summation is a $q$-analogue of \eqref{2f1}
(cf.\ \cite[Appendix~(II.6)]{GR}).
\begin{equation}\label{2f1q}
{}_2\phi_1\!\left[\begin{matrix}a,q^{-n}\\c\end{matrix};q,q\right]=
\frac{(c/a;q)_n}{(c;q)_n}\,a^n.
\end{equation}

We provide the following noncommutative extensions of \eqref{2f1q}:
\begin{Theorem}\label{nc2f1Q}
Let $A$ and $C$ be noncommutative parameters of some unit ring,
and suppose that $Q$ commutes both with $A$ and $C$. Further,
let $n$ be a nonnegative integer. Then we have the following
summation for a noncommutative hypergeometric series of type I.
\begin{equation}\label{nc2f1QIgl}
{}_2\phi_1\!\left\lceil\begin{matrix}A,Q^{-n}\\
C\end{matrix};Q,Q\right\rfloor=
\left\lceil\begin{matrix}CA^{-1}\\C\end{matrix};Q,A\right\rfloor_n.
\end{equation}
Further, we we have the following
summation for a noncommutative hypergeometric series of type II.
\begin{equation}\label{nc2f1QIIgl}
{}_2\phi_1\!\left\lfloor\begin{matrix}A,Q^{-n}\\
C\end{matrix};Q,Q\right\rceil=
\left\lfloor\begin{matrix}CA^{-1}Q\\C\end{matrix};Q,A\right\rceil_n
(A-CQ^n)^{-1}(A-C).
\end{equation}
\end{Theorem}

Note that the right-hand side of \eqref{nc2f1QIIgl}
may also be written as
\begin{equation*}
\left\lfloor\begin{matrix}CA^{-1}Q\\C\end{matrix};Q,A\right\rceil_{n-1}
(I-CQ^{n-1})^{-1}(A-C),
\end{equation*}
however, as in \eqref{nc2f1IIgl} it is in the simple case $n=0$
easier to see that the sum reduces correctly to $I$ when one
uses the expression on the right-hand side of \eqref{nc2f1QIIgl}.

\begin{proof}[Proof of Theorem~\ref{nc2f1Q}]
We prove \eqref{nc2f1QIgl} by induction on $n$ (and leave
\eqref{nc2f1QIIgl} to the reader).
For $n=0$ the formula is trivial.
Assume the formula is true up to a fixed $n$.
To prove it for $n+1$, use the elementary identity
\begin{equation*}
\left\lceil\begin{matrix}Q^{-(n+1)}\\Q\end{matrix};Q,I\right\rfloor_k=
\left\lceil\begin{matrix}Q^{-n}\\Q\end{matrix};Q,I\right\rfloor_k
\left[I-Q^{-n-1}(I-Q^{-n-1+k})^{-1}(I-Q^k)\right]
\end{equation*}
to obtain
\begin{multline*}
{}_2\phi_1\!\left\lceil\begin{matrix}A,Q^{-(n+1)}\\
C\end{matrix};Q,Q\right\rfloor\\=
\sum_{k=0}^n
\left\lceil\begin{matrix}A,Q^{-n}\\C,Q\end{matrix};Q,Q\right\rfloor_k
\left[I-Q^{-n-1}(I-Q^{-n-1+k})^{-1}(I-Q^k)\right]\\=
{}_2\phi_1\!\left\lceil\begin{matrix}A,Q^{-n}\\
C\end{matrix};Q,Q\right\rfloor-
{}_2\phi_1\!\left\lceil
\begin{matrix}AQ,Q^{-n}\\CQ\end{matrix};Q,Q\right\rfloor
Q^{-n}(I-C)^{-1}(I-A)\\=
\left\lceil\begin{matrix}CA^{-1}\\C\end{matrix};Q,A\right\rfloor_n-
\left\lceil\begin{matrix}CA^{-1}\\
CQ\end{matrix};Q,A\right\rfloor_n(I-C)^{-1}(I-A),
\end{multline*}
the last equation due to the inductive hypothesis.
Now apply Lemma~\ref{lem21}, Equation~\eqref{lem21QIgl},
which completes the proof of \eqref{nc2f1QIgl}.
\end{proof}

Here is another $q$-analogue of \eqref{2f1} (cf.\ \cite[Appendix~(II.7)]{GR}).
\begin{equation}\label{2f1qr}
{}_2\phi_1\!\left[\begin{matrix}a,q^{-n}\\
c\end{matrix};q,\frac{cq^n}a\right]=
\frac{(c/a;q)_n}{(c;q)_n}.
\end{equation}

We provide the following noncommutative extensions of \eqref{2f1qr}:
\begin{Theorem}\label{nc2f1Qr}
Let $A$ and $C$ be noncommutative parameters of some unit ring,
and suppose that $Q$ commutes both with $A$ and $C$. Further,
let $n$ be a nonnegative integer. Then we have the following
summation for a noncommutative hypergeometric series of type I.
\begin{equation}\label{nc2f1QIrgl}
{}_2\phi_1\!\left\lceil\begin{matrix}A,Q^{-n}\\
C\end{matrix};Q,A^{-1}CQ^n\right\rfloor=
\left\lceil\begin{matrix}A^{-1}C\\C\end{matrix};Q,I\right\rfloor_n.
\end{equation}
Further, we we have the following
summation for a noncommutative hypergeometric series of type II.
\begin{equation}\label{nc2f1QIIrgl}
{}_2\phi_1\!\left\lfloor\begin{matrix}A,Q^{-n}\\
C\end{matrix};Q,A^{-1}CQ^n\right\rceil=
\left\lfloor\begin{matrix}A^{-1}CQ\\C\end{matrix};Q,I\right\rceil_n
(I-A^{-1}CQ^n)^{-1}(I-A^{-1}C).
\end{equation}
\end{Theorem}

\begin{proof}
We prove \eqref{nc2f1QIrgl} by induction on $n$ (and leave
\eqref{nc2f1QIIrgl} to the reader).
For $n=0$ the formula is trivial.
Assume the formula is true up to a fixed $n$.
To prove it for $n+1$, use the elementary identity
\begin{equation*}
\left\lceil\begin{matrix}Q^{-(n+1)}\\Q\end{matrix};Q,Q^{n+1}\right\rfloor_k=
\left\lceil\begin{matrix}Q^{-n}\\Q\end{matrix};Q,Q^n\right\rfloor_k
\left[I-(I-Q^{-n-1+k})^{-1}(I-Q^k)\right]
\end{equation*}
to obtain
\begin{multline*}
{}_2\phi_1\!\left\lceil\begin{matrix}A,Q^{-(n+1)}\\
C\end{matrix};Q,A^{-1}CQ^{n+1}\right\rfloor=
\sum_{k=0}^n
\left\lceil\begin{matrix}A,Q^{-n}\\
C,Q\end{matrix};Q,A^{-1}CQ^n\right\rfloor_k\\\times
\left[I-(I-Q^{-n-1+k})^{-1}(I-Q^k)\right]=
{}_2\phi_1\!\left\lceil\begin{matrix}A,Q^{-n}\\
C\end{matrix};Q,A^{-1}CQ^n\right\rfloor\\-
{}_2\phi_1\!\left\lceil\begin{matrix}AQ,Q^{-n}\\
CQ\end{matrix};Q,A^{-1}CQ^n\right\rfloor
(I-C)^{-1}(I-A)A^{-1}CQ^n\\=
\left\lceil\begin{matrix}A^{-1}C\\C\end{matrix};Q,I\right\rfloor_n-
\left\lceil\begin{matrix}A^{-1}C\\
CQ\end{matrix};Q,I\right\rfloor_n(I-C)^{-1}(I-A)A^{-1}CQ^n,
\end{multline*}
the last equation due to the inductive hypothesis.
Now apply Lemma~\ref{lem21}, Equation~\eqref{lem21QIrgl},
and the proof of \eqref{nc2f1QIrgl} is complete.
\end{proof}

There are two ways to obtain Theorem~\ref{nc2f1Qr} from
Theorem~\ref{nc2f1Q} directly, namely by {\em inverting the basis}
$Q\to Q^{-1}$, or by {\em reversing the sum}. We give some
details for both of these possibilities, but restrict ourserves
to the derivation of \eqref{nc2f1QIrgl} from \eqref{nc2f1QIgl}.
The details of deriving \eqref{nc2f1QIIrgl} from \eqref{nc2f1QIIgl}
by inverting the basis or by reversing the sum are similar
and left left to the reader.

\smallskip
\noindent 1. {\em Inverting the base.} It is easy to verify that
\begin{equation}\label{invbase}
\left\lceil\begin{matrix}A\\C\end{matrix};Q^{-1},Z\right\rfloor_k=
C^{-1}\left\lceil\begin{matrix}A^{-1}\\
C^{-1}\end{matrix};Q,AZC^{-1}\right\rfloor_k C.
\end{equation}
Now since
\begin{equation*}
{}_2\phi_1\!\left\lceil\begin{matrix}A,Q^n\\
C\end{matrix};Q^{-1},Q^{-1}\right\rfloor=
\left\lceil\begin{matrix}CA^{-1}\\C\end{matrix};Q^{-1},A\right\rfloor_n
\end{equation*}
by \eqref{nc2f1QIgl}, we deduce from \eqref{invbase} that
\begin{equation*}
C^{-1}\,{}_2\phi_1\!\left\lceil\begin{matrix}A^{-1},Q^{-n}\\
C^{-1}\end{matrix};Q,AC^{-1}Q^n\right\rfloor C=
C^{-1}\left\lceil\begin{matrix}AC^{-1}\\
C^{-1}\end{matrix};Q,I\right\rfloor_n C.
\end{equation*}
Now simply replace $A$ by $A^{-1}$ and $C$ by $C^{-1}$ and then multiply
both sides by $C^{-1}$ from the left and by $C$ from the right to
get \eqref{nc2f1QIrgl}.

\smallskip
\noindent 2. {\em Reversing the sum.}
Using
\begin{align*}
\prod_{j=1}^{n-k}B_{n-k-j}=\prod_{j=1}^{-k}B_{n-k-j}
\prod_{j=1-k}^{n-k}B_{n-k-j}&{}=\prod_{j=1-k}^0B_{n-1+j}^{-1}
\prod_{j=1}^nB_{n-j}\\
&{}=\prod_{j=1}^kB_{n-1-k+j}^{-1}\prod_{j=1}^nB_{n-j}
\end{align*}
which is derived using \eqref{invprod}, we readily deduce
\begin{equation*}
\left\lceil\begin{matrix}A\\C\end{matrix};Q,I\right\rfloor_{n-k}=
A^{-1}\left\lceil\begin{matrix}C^{-1}Q^{1-n}\\
A^{-1}Q^{1-n}\\\end{matrix};Q,CA^{-1}\right\rfloor_k A
\left\lceil\begin{matrix}A\\C\end{matrix};Q,I\right\rfloor_n.
\end{equation*}
Similarly,
\begin{align*}
\left\lceil\begin{matrix}Q^{-n}\\Q\end{matrix};Q,I\right\rfloor_{n-k}
&{}=\left\lceil\begin{matrix}Q^{-n}\\
Q\end{matrix};Q,I\right\rfloor_k
\left\lceil\begin{matrix}Q^{-n}\\Q\end{matrix};Q,I\right\rfloor_n
Q^{(n+1)k}\\&{}=
\left\lceil\begin{matrix}Q^{-n}\\
Q\end{matrix};Q,I\right\rfloor_k (-1)^n\,Q^{(n+1)k-\binom{n+1}2}.
\end{align*}
Hence we have from \eqref{nc2f1QIgl}, by reversing the sum
on the left hand side,
\begin{align*}
A^{-1}{}_2\phi_1\!\left\lceil\begin{matrix}C^{-1}Q^{1-n},Q^{-n}\\
A^{-1}Q^{1-n}\end{matrix};Q,CA^{-1}Q^n\right\rfloor
A\,(-1)^n\,Q^{-\binom n2}\left\lceil\begin{matrix}A\\
C\end{matrix};Q,I\right\rfloor_n&\\=
\left\lceil\begin{matrix}CA^{-1}\\C\end{matrix};Q,A\right\rfloor_n&.
\end{align*}
Performing the simultaneous substitutions $A\mapsto C^{-1}Q^{1-n}$,
$C\mapsto A^{-1}Q^{1-n}$, and putting some factors to the other side
gives
\begin{multline}\label{aux1}
{}_2\phi_1\!\left\lceil\begin{matrix}A,Q^{-n}\\
C\end{matrix};Q,A^{-1}CQ^n\right\rfloor\\=
C^{-1}\left\lceil\begin{matrix}A^{-1}C\\
A^{-1}Q^{1-n}\end{matrix};Q,C^{-1}Q^{1-n}\right\rfloor_n
\left\lceil\begin{matrix}C^{-1}Q^{1-n}\\
A^{-1}Q^{1-n}\end{matrix};Q,I\right\rfloor_n^{-1}
C\,(-1)^n\,Q^{\binom n2}.
\end{multline}
We want to show that the right-hand side of \eqref{aux1} reduces to
the right-hand side of \eqref{nc2f1QIrgl}.
First, we compute
\begin{equation}\label{aux3}
\left\lceil\begin{matrix}C^{-1}Q^{1-n}\\
A^{-1}Q^{1-n}\end{matrix};Q,I\right\rfloor_n^{-1}=
C\left\lceil\begin{matrix}A\\C\end{matrix};Q,A^{-1}C\right\rfloor_n
C^{-1}.
\end{equation}
Further,
\begin{multline}\label{aux4}
\left\lceil\begin{matrix}A^{-1}C\\
A^{-1}Q^{1-n}\end{matrix};Q,C^{-1}Q^{1-n}\right\rfloor_n\\=
(-1)^n\,Q^{\binom n2}\prod_{j=1}^n\left[
A(I-AQ^{j-1})^{-1}(I-A^{-1}CQ^{n-j})C^{-1}\right].
\end{multline}
Thus, applying \eqref{aux3} and \eqref{aux4} to
the right-hand side of \eqref{aux1} and equating the result to
the right-hand side of \eqref{nc2f1QIrgl}, we have established
\eqref{nc2f1QIrgl} once we have shown the following lemma.

\begin{Lemma}\label{lem1}
\begin{equation}\label{lem1gl}
\prod_{j=1}^n\left[C^{-1}
A(I-AQ^{j-1})^{-1}(I-A^{-1}CQ^{n-j})\right]
\left\lceil\begin{matrix}A\\C\end{matrix};Q,A^{-1}C\right\rfloor_n
=\left\lceil\begin{matrix}A^{-1}C\\C\end{matrix};Q,I\right\rfloor_n.
\end{equation}
\end{Lemma}

\begin{proof}
We proceed by induction on $n$. For $n=0$ the identity is trivial.
Observe that for $n=1$ the statement amounts to
\begin{equation}\label{aux2}
C^{-1}A(I-A)^{-1}(I-A^{-1}C)(I-C)^{-1}(I-A)A^{-1}C=
(I-C)^{-1}(I-A^{-1}C).
\end{equation}
This is easily verified by splitting the factor
$(I-A^{-1}C)(I-C)^{-1}$ on the left-hand side of \eqref{aux2}
in two terms as $I-(I-A)A^{-1}C(I-C)^{-1}$, and simplifying
the expression to $I-(I-C)^{-1}(I-A)A^{-1}C$
which clearly equals $(I-C)^{-1}(I-A^{-1}C)$.

Assume now that \eqref{lem1gl} is true for all nonnegative integers
up to a fixed $n$. Then
\begin{multline}
\prod_{j=1}^n\left[C^{-1}
A(I-AQ^{j-1})^{-1}(I-A^{-1}CQ^{n-j})\right]
\left\lceil\begin{matrix}A\\C\end{matrix};Q,A^{-1}C\right\rfloor_n\\=
C^{-1}A(I-A)^{-1}(I-A^{-1}CQ^{n-1})
\prod_{j=1}^{n-1}\left[C^{-1}
A(I-AQ^j)^{-1}(I-A^{-1}CQ^{n-1-j})\right]\\\times
\left\lceil\begin{matrix}AQ\\CQ\end{matrix};Q,A^{-1}C\right\rfloor_{n-1}
(I-C)^{-1}(I-A)A^{-1}C\\=
C^{-1}A(I-A)^{-1}(I-A^{-1}CQ^{n-1})
\left\lceil\begin{matrix}A^{-1}C\\CQ\end{matrix};Q,I\right\rfloor_{n-1}
(I-C)^{-1}(I-A)A^{-1}C,
\end{multline}
the last equation due to the $A\mapsto AQ$, $C\mapsto CQ$ case of
the inductive hypothesis. The proof is complete after application
of Lemma~\ref{lem2}.
\end{proof}

\begin{Lemma}\label{lem2}
\begin{multline}\label{lem2gl}
C^{-1}A(I-A)^{-1}(I-A^{-1}CQ^{n-1})
\left\lceil\begin{matrix}A^{-1}C\\CQ\end{matrix};Q,I\right\rfloor_{n-1}
(I-C)^{-1}(I-A)A^{-1}C\\
=\left\lceil\begin{matrix}A^{-1}C\\C\end{matrix};Q,I\right\rfloor_n.
\end{multline}
\end{Lemma}

\begin{proof}
We proceed by induction on $n$. For $n=1$ \eqref{lem2gl} is simply
\eqref{aux2}.
Assume now that \eqref{lem2gl} is true for all positive integers
up to a fixed $n$. Then
\begin{multline}
C^{-1}A(I-A)^{-1}(I-A^{-1}CQ^n)
\left\lceil\begin{matrix}A^{-1}C\\CQ\end{matrix};Q,I\right\rfloor_n
(I-C)^{-1}(I-A)A^{-1}C\\=
C^{-1}Q^{-1}A(I-A)^{-1}(I-A^{-1}CQ^n)
\left\lceil\begin{matrix}A^{-1}CQ\\CQ^2\end{matrix};
Q,I\right\rfloor_{n-1}\\\times
(I-CQ)^{-1}(I-A^{-1}C)(I-C)^{-1}(I-A)A^{-1}CQ\\=
\left\lceil\begin{matrix}A^{-1}CQ\\CQ\end{matrix};Q,I\right\rfloor_{n-1}
C^{-1}Q^{-1}A(I-A)^{-1}(I-A^{-1}C)(I-C)^{-1}(I-A)A^{-1}CQ,
\end{multline}
the last equation
due to the $C\mapsto CQ$ case of the inductive hypothesis. Finally,
an application of \eqref{aux2} simplifies the last expression to
\begin{equation*}
\left\lceil\begin{matrix}A^{-1}CQ\\CQ\end{matrix};Q,I\right\rfloor_{n-1}
(I-C)^{-1}(I-A^{-1}C)
\end{equation*}
which is equal to the generalized $Q$-shifted factorial on the
right-hand side of \eqref{lem2gl}.
\end{proof}

\section{Pfaff--Saalsch\"utz summations}\label{sec3f2}

The classical Pfaff--Saalsch\"utz summation formula
(cf.\ \cite[Appendix~(III.2)]{S})
sums a terminating balanced ${}_3F_2$ series:
\begin{equation}\label{3f2}
{}_3F_2\!\left[\begin{matrix}a,b,-n\\c,a+b-c+1-n\end{matrix};1\right]=
\frac{(c-a)_n(c-b)_n}{(c)_n(c-a-b)_n}.
\end{equation}

In order to derive a noncommutative extension of \eqref{3f2},
one should at least be able to extend its $n=1$ special case,
which is
\begin{equation*}
1-\frac{ab}{c(a+b-c)}=\frac{(c-a)(c-b)}{c(c-a-b)}.
\end{equation*}

This is no problem indeed, as we have
\begin{equation}\label{easy32}
I-C^{-1}A(A+B-C)^{-1}B=C^{-1}(C-B)(C-A-B)^{-1}(C-A),
\end{equation}
for noncommutative parameters $A$, $B$, $C$,
as one immediately verifies. In fact,
\begin{multline*}
I-C^{-1}A(A+B-C)^{-1}B=C^{-1}[C-A(A+B-C)^{-1}B]\\
=C^{-1}[C-((A+B-C)-(B-C))(A+B-C)^{-1}B]\\
=C^{-1}[C-B-(C-B)(A+B-C)^{-1}B]\\
=C^{-1}(C-B)[I-(A+B-C)^{-1}B]\\=
C^{-1}(C-B)(C-A-B)^{-1}(C-A).
\end{multline*}

Nevertheless, we were not able to extend \eqref{easy32}
to a noncommutative Pfaff--Saalsch\"utz summation with
arbitrary noncommutative parameters $A$, $B$ and $C$.
For the case of general $n$, we need a restriction on the
sum $D=A+B-C$, namely that it must commute with the other elements
$A$, $B$ and $C$ (e.g., $D=dI$).

Here are two noncommutative extensions of \eqref{3f2}:
\begin{Theorem}\label{nc3f2}
Let $A$, $B$ and $C$ be noncommutative parameters of some unit ring,
and assume that the sum $A+B-C$ commutes each with $A$, $B$ and $C$.
Further, let $n$ be a nonnegative integer. Then we have the following
summation for a noncommutative hypergeometric series of type I.
\begin{equation}\label{nc3f2Igl}
{}_3F_2\!\left\lceil\begin{matrix}A,B,-nI\\
C,A+B-C+(1-n)I\end{matrix}\,;I\right\rfloor=
\left\lceil\begin{matrix}C-B,C-A\\
C,C-A-B\end{matrix}\right\rfloor_n.
\end{equation}
Further, we have the following
summation for a noncommutative hypergeometric series of type II.
\begin{multline}\label{nc3f2IIgl}
{}_3F_2\!\left\lfloor\begin{matrix}A,B,-nI\\
C,A+B-C+(1-n)I\end{matrix}\,;I\right\rceil=
\left\lfloor\begin{matrix}C-B+I,C-A+I\\
C,C-A-B\end{matrix}\right\rceil_n\\\times
(C-A+nI)^{-1}(C-B+nI)^{-1}(C-B)(C-A).
\end{multline}
\end{Theorem}
Note that the right-hand side of \eqref{nc3f2IIgl} may also be written as
\begin{equation*}
\left\lfloor\begin{matrix}C-B+I,C-A+I\\C,C-A-B\end{matrix}\right\rceil_{n-1}
(C+(n-1)I)^{-1}(C-A-B+(n-1)I)^{-1}(C-B)(C-A).
\end{equation*}

\begin{proof}[Proof of Theorem~\ref{nc3f2}]
We prove \eqref{nc3f2IIgl} by induction on $n$,
leaving the proof of \eqref{nc3f2Igl} (which is similar) to the reader.
For $n=0$ the formula is trivial.
Assume the formula is true up to a fixed $n$.
To prove it for $n+1$, use the elementary identity
\begin{multline*}
\left\lfloor\begin{matrix}-(n+1)I\\A+B-C-nI\end{matrix}\right\rceil_k=
\left\lfloor\begin{matrix}-nI\\A+B-C+(1-n)I\end{matrix}\right\rceil_k\\\times
\left[I-(A+B-C-nI)^{-1}(A+B-C+I)((-n-1+k)I)^{-1}kI\right]
\end{multline*}
to obtain
\begin{multline*}
{}_3F_2\!\left\lfloor\begin{matrix}A,B,-(n+1)I\\
C,A+B-C-nI\end{matrix}\,;I\right\rceil=
\sum_{k=0}^n\left\lfloor\begin{matrix}A,B,-nI\\
C,A+B-C+(1-n)I,I\end{matrix}\right\rceil_k\\\times
\left[I-(A+B-C-nI)^{-1}(A+B-C+I)((-n-1+k)I)^{-1}kI\right]\\=
{}_3F_2\!\left\lfloor\begin{matrix}A,B,-nI\\
C,A+B-C+(1-n)I\end{matrix}\,;I\right\rceil\\-
C^{-1}A(A+B-C+(1-n)I)^{-1}B(A+B-C-nI)^{-1}(A+B-C+I)\\\times
{}_3F_2\!\left\lfloor\begin{matrix}A+I,B+I,-nI\\
C+I,A+B-C+(2-n)I\end{matrix}\,;I\right\rceil\\=
\left\lfloor\begin{matrix}C-B+I,C-A+I\\
C,C-A-B\end{matrix}\right\rceil_n
(C-A+nI)^{-1}(C-B+nI)^{-1}(C-B)(C-A)\\
-C^{-1}A(A+B-C+(1-n)I)^{-1}B(A+B-C-nI)^{-1}(A+B-C+I)\\\times
\left\lfloor\begin{matrix}C-B+I,C-A+I\\
C+I,C-A-B-I\end{matrix}\right\rceil_n
(C-A+nI)^{-1}(C-B+nI)^{-1}(C-B)(C-A),
\end{multline*}
the last equation due to the inductive hypothesis.
We are done after application of Lemma~\ref{lem32},
Equation~\eqref{lem32IIgl}.
\end{proof}

The $q$-analogue of \eqref{3f2} is the following summation formula
(cf.\ \cite[Appendix~(II.12)]{GR}).
\begin{equation}\label{3f2q}
{}_3\phi_2\!\left[\begin{matrix}a,b,q^{-n}\\
c,abq^{1-n}/c\end{matrix};q,q\right]=
\frac{(c/a;q)_n(c/b;q)_n}{(c;q)_n(c/ab;q)_n}.
\end{equation}

In order to derive a noncommutative extension of \eqref{3f2q},
one should at least be able to extend its $n=1$ special case,
which is
\begin{equation*}
1-\frac{(1-a)(1-b)}{(1-c)(1-ab/c)}=\frac{(1-c/a)(1-c/b)}{(1-c)(1-c/ab)}.
\end{equation*}

This is no problem indeed, as we have
\begin{multline}\label{easy32Q}
I-(I-C)^{-1}(I-A)(I-BC^{-1}A)^{-1}(I-B)\\=
(I-C)^{-1}(I-CB^{-1})(I-A^{-1}CB^{-1})^{-1}(I-A^{-1}C),
\end{multline}
for noncommutative parameters $A$, $B$, $C$,
as one immediately verifies. In fact,
\begin{multline*}
I-(I-C)^{-1}(I-A)(I-BC^{-1}A)^{-1}(I-B)\\
=(I-C)^{-1}[I-C-(I-A)(I-BC^{-1}A)^{-1}(I-B)]\\
=(I-C)^{-1}[I-C-((I-BC^{-1}A)-(I-BC^{-1})A)(I-BC^{-1}A)^{-1}(I-B)]\\
=(I-C)^{-1}[I-C-(I-B)-(I-CB^{-1})BC^{-1}A(I-BC^{-1}A)^{-1}(I-B)]\\
=(I-C)^{-1}(I-CB^{-1})[B+(I-A^{-1}CB^{-1})^{-1}(I-B)]\\
=(I-C)^{-1}(I-CB^{-1})(I-A^{-1}CB^{-1})^{-1}[B-A^{-1}C+(I-B)]\\
=(I-C)^{-1}(I-CB^{-1})(I-A^{-1}CB^{-1})^{-1}(I-A^{-1}C).
\end{multline*}

Nevertheless, we were not able to extend \eqref{easy32Q}
to a noncommutative $q$-Pfaff--Saalsch\"utz summation with
arbitrary noncommutative parameters $A$, $B$ and $C$.
For the case of general $n$, we need a restriction on the
product $D=BC^{-1}A$, namely that it must commute with the other elements
$A$, $B$ and $C$ (e.g., $D=dI$).

Here are two noncommutative extensions of \eqref{3f2q}:
\begin{Theorem}\label{nc3f2Q}
Let $A$, $B$ and $C$ be noncommutative parameters of some unit ring,
and suppose that $Q$ commutes each with $A$, $B$ and $C$. Moreover, assume
that the product $BC^{-1}A$ commutes each with $A$, $B$ and $C$.
Further, let $n$ be a nonnegative integer. Then we have the following
summation for a noncommutative basic hypergeometric series of type I.
\begin{equation}\label{nc3f2QIgl}
{}_3\phi_2\!\left\lceil\begin{matrix}A,B,Q^{-n}\\
C,BC^{-1}AQ^{1-n}\end{matrix}\,;Q,Q\right\rfloor=
\left\lceil\begin{matrix}CB^{-1},A^{-1}C\\
C,A^{-1}CB^{-1}\end{matrix};Q,I\right\rfloor_n.
\end{equation}
Further, we have the following
summation for a noncommutative basic hypergeometric series of type II.
\begin{multline}\label{nc3f2QIIgl}
{}_3\phi_2\!\left\lfloor\begin{matrix}A,B,Q^{-n}\\
C,BC^{-1}AQ^{1-n}\end{matrix}\,;Q,Q\right\rceil=
\left\lfloor\begin{matrix}CB^{-1}Q,A^{-1}CQ\\
C,A^{-1}CB^{-1}\end{matrix};Q,I\right\rceil_n\\\times
(I-A^{-1}CQ^n)^{-1}(I-CB^{-1}Q^n)^{-1}(I-CB^{-1})(I-A^{-1}C).
\end{multline}
\end{Theorem}
Note that the right-hand side of \eqref{nc3f2QIIgl} may also be written as
\begin{multline*}
\left\lfloor\begin{matrix}CB^{-1}Q,A^{-1}CQ\\
C,A^{-1}CB^{-1}\end{matrix};Q,I\right\rceil_{n-1}\\\times
(I-CQ^{n-1})^{-1}(I-A^{-1}CB^{-1}Q^{n-1})^{-1}(I-CB^{-1})(I-A^{-1}C).
\end{multline*}

\begin{proof}[Proof of Theorem~\ref{nc3f2Q}]
We prove \eqref{nc3f2QIIgl} by induction on $n$.
We leave the proof of \eqref{nc3f2QIgl} (which is similar) to the reader.
For $n=0$ the formula is trivial.
Assume the formula is true up to a fixed $n$.
To prove it for $n+1$, use the elementary identity
\begin{multline*}
\left\lfloor\begin{matrix}Q^{-(n+1)}\\
BC^{-1}AQ^{-n}\end{matrix};Q,I\right\rceil_k=
\left\lfloor\begin{matrix}Q^{-n}\\
BC^{-1}AQ^{1-n}\end{matrix};Q,I\right\rceil_k\\\times
\left[I-(I-BC^{-1}AQ^{-n})^{-1}(I-BC^{-1}AQ)
(I-Q^{-n-1+k})^{-1}(I-Q^k)\right]
\end{multline*}
to obtain
\begin{multline*}
{}_3\phi_2\!\left\lfloor\begin{matrix}A,B,Q^{-(n+1)}\\
C,BC^{-1}AQ^{-n}\end{matrix}\,;Q,Q\right\rceil=
\sum_{k=0}^n\left\lfloor\begin{matrix}A,B,Q^{-n}\\
C,BC^{-1}AQ^{1-n},Q\end{matrix}\,;Q,Q\right\rceil_k\\\times
\left[I-(I-BC^{-1}AQ^{-n})^{-1}(I-BC^{-1}AQ)
(I-Q^{-n-1+k})^{-1}(I-Q^k)\right]\\=
{}_3\phi_2\!\left\lfloor\begin{matrix}A,B,Q^{-n}\\
C,BC^{-1}AQ^{1-n}\end{matrix}\,;Q,Q\right\rceil\\-
(I-C)^{-1}(I-A)(I-BC^{-1}AQ^{1-n})^{-1}(I-B)
(I-BC^{-1}AQ^{-n})^{-1}(I-BC^{-1}AQ)\\\times
{}_3\phi_2\!\left\lfloor\begin{matrix}AQ,BQ,Q^{-n}\\
CQ,BC^{-1}AQ^{2-n}\end{matrix}\,;Q,Q\right\rceil\\=
\left\lfloor\begin{matrix}CB^{-1}Q,A^{-1}CQ\\
C,A^{-1}CB^{-1}\end{matrix};Q,I\right\rceil_n\\\times
(I-A^{-1}CQ^n)^{-1}(I-CB^{-1}Q^n)^{-1}(I-CB^{-1})(I-A^{-1}C)\\
-(I-C)^{-1}(I-A)(I-BC^{-1}AQ^{1-n})^{-1}(I-B)
(I-BC^{-1}AQ^{-n})^{-1}(I-BC^{-1}AQ)\\\times
\left\lfloor\begin{matrix}CB^{-1}Q,A^{-1}CQ\\
CQ,A^{-1}CB^{-1}Q^{-1}\end{matrix};Q,I\right\rceil_n\\\times
(I-A^{-1}CQ^n)^{-1}(I-CB^{-1}Q^n)^{-1}(I-CB^{-1})(I-A^{-1}C),
\end{multline*}
the last equation due to the inductive hypothesis.
We are done after application of Lemma~\ref{lem32Q},
Equation~\ref{lem32QIIgl}.
\end{proof}

\section{Very-well-poised ${}_7F_6$ summations}\label{sec7f6}

Dougall's summation formula (cf.\ \cite[Appendix~(III.14)]{S})
sums a terminating very-well-poised $2$-balanced ${}_7F_6$ series:
\begin{multline}\label{7f6}
{}_7F_6\!\left[\begin{matrix}\frac a2+1,a,b,c,d,2a-b-c-d+1+n,-n\\
\frac a2,a-b+1,a-c+1,a-d+1,b+c+d-a-n,a+n+1\end{matrix};1\right]\\=
\frac{(a+1)_n(a-b-c+1)_n(a-b-d+1)_n(a-c-d+1)_n}
{(a-b+1)_n(a-c+1)_n(a-d+1)_n(a-b-c-d+1)_n}.
\end{multline}

In order to derive a noncommutative extension of \eqref{7f6},
one should at least be able to extend its $n=1$ special case,
which (with $a$ replaced by $a-1$) is
\begin{multline}\label{7f61}
1-\frac{bcd(2a-b-c-d)}{(a-b)(a-c)(a-d)(b+c+d-a)}\\=
\frac{a(a-b-c)(a-b-d)(a-c-d)}{(a-b)(a-c)(a-d)(a-b-c-d)}.
\end{multline}

A noncommutative extension of \eqref{7f61} is indeed available
when one assumes that $A$, $B$, $C$, and $D$ are parameters such that
$A$, $B$ and $D$ mutually commute while $C$ does not commute with
any of the other parameters. We then have
\begin{multline}\label{easy76}
I-(A-C)^{-1}B(A-B)^{-1}C(B+C+D-A)^{-1}D(A-D)^{-1}(2A-B-C-D)\\
=(A-C)^{-1}(A-B-D)(A-B)^{-1}(A-C-D)\\\times
(A-B-C-D)^{-1}A(A-D)^{-1}(A-B-C).
\end{multline}
We already know from \eqref{easy32} that
\begin{equation}\label{easy32G}
I-G^{-1}E(E+F-G)^{-1}F=G^{-1}(G-F)(G-E-F)^{-1}(G-E)
\end{equation}
holds, for noncommuting parameters $E$, $F$ and $G$.
Now let
\begin{align*}
E&= BC,\\
F&= D(2A-B-C-D),\\
G&= (A-B)(A-C),
\end{align*}
while assuming that
that $A$, $B$, $C$, and $D$ are parameters such that
$A$, $B$ and $D$ mutually commute while $C$ does not commute with
any of the other parameters.
Then one readily verifies that
\begin{align*}
G-E&= A(A-B-C),\\
G-F&= (A-B-D)(A-C-D),\\
G-E-F&= (A-D)(A-B-C-D).
\end{align*}
After these substitutions into \eqref{easy32G}
and rearranging the order of some factors (which do not involve $C$),
we immediately establish \eqref{easy76}.

However, rather than \eqref{easy76}, we prefer to consider its variant
which is obtained by replacing $D$ by $2A-B-C-D$, for convenience.
The result is the following.

Under the assumption that $A$, $B$, $C$, and $D$ are (noncommuting)
parameters such that $A$, $B$ and $C+D$ mutually commute, there holds
\begin{multline}\label{easy76v}
I-(A-C)^{-1}B(A-B)^{-1}C(A-D)^{-1}(2A-B-C-D)(B+C+D-A)^{-1}D\\
=(A-C)^{-1}(A-C-D)(A-B)^{-1}(A-B-D)\\\times
(A-D)^{-1}A(A-B-C-D)^{-1}(A-B-C).
\end{multline}

As a matter of fact, we were not able to extend \eqref{easy76}
(or the equivalent \eqref{easy76v}) to a noncommutative extension
of Dougall's summation \eqref{7f6}, valid for general $n$,
without introducing more commutation relations.
What we require in the general case (in place of the weaker
requirement that $A$, $B$ and $C+D$ shall mutually commute)
are the following commutation relations,
\begin{subequations}\label{comrel}
\begin{align}\label{comrela}
A\quad &\text{commutes with \quad$B$, $C$, $D$,}\\\label{comrelbcdp}
B+C+D\quad &\text{commutes with \quad$B$, $C$, $D$,}
\end{align}
and the following ``rotation relations'',
\begin{equation}\label{comrelrot}
BCD=CDB=DBC.
\end{equation}
\end{subequations}

It is easy to check that
\eqref{comrelbcdp} implies
\begin{equation}\label{comrelbcd}
BC-CB=CD-DC=DB-BD,
\end{equation}
which, using Lie brackets, writes elegantly as
\begin{equation*}
[B,C]=[C,D]=[D,B].
\end{equation*}
Of course, these Lie products need not to be $O$
(the zero element of the unit ring $R$).
Suppose that $E$ commutes each with $B$, $C$ and $D$.
Since it follows from \eqref{comrelrot} that
\begin{equation*}
BCD\quad \text{commutes with}\quad BC+BD+CD,
\end{equation*}
we immediately deduce from
\begin{equation*}
(B+E)(C+E)(D+E)=BCD+(BC+BD+CD)E+(B+C+D)E^2+E^3,
\end{equation*}
that \eqref{comrelbcdp} and \eqref{comrelrot} imply that
\begin{equation}\label{comrelbcdd}
BCD\quad \text{commutes with \quad $(B+E)(C+E)(D+E)$.}
\end{equation}

\begin{Remark}\label{rem}
One may wonder whether the above commutation relations \eqref{comrel}
allow any room for noncommutativity, where $B$, $C$, $D$ do not already
mutually commute. The following nontrivial example of $A$, $B$, $C$, $D$
realizing \eqref{comrel} has been kindly communicated to us
(essentially in the given form) by Hjalmar Rosengren.

Consider $R$ to be the ring of $(3\times 3)$-matrices over the complex
numbers $\mathbb C$ with unit element $I$ and zero element $O$.
Let $\omega=e^{2\pi i/3}$, a cubic root of unity.
Suppose $a,b,c,d\in\mathbb C$ with $b\neq 0$, $c\neq 0$. Let $A=aI$,
\begin{equation*}
B=\begin{pmatrix}d&b&0\\b^{-1}&d&c\\0&-\omega^2c^{-1}&d\end{pmatrix},
C=\begin{pmatrix}d&\omega b&0\\
\omega^2 b^{-1}&d&\omega c\\0&-\omega c^{-1}&d\end{pmatrix},
D=\begin{pmatrix}d&\omega^2 b&0\\
\omega b^{-1}&d&\omega^2 c\\0&-c^{-1}&d\end{pmatrix}.
\end{equation*}
Then the matrices $A$, $B$, $C$ and $D$ satisfy \eqref{comrel}:
Of course, \eqref{comrela} is trivially satisfied.
Observe that $B$, $C$, $D$ do not mutually commute since
\begin{equation*}
X=BC-CB=CD-DC=DB-BD=(\omega^2-\omega)\begin{pmatrix}1&0&0\\0&\omega&0\\
0&0&\omega^2\end{pmatrix}\neq O.
\end{equation*}
Further observe that with the above choice of $B$, $C$ and $D$
one has $C=X^{-1}BX$, $D=XBX^{-1}$, hence $XC=BX$, $DX=XB$, and $CX=XD$.
It is straightforward to compute $B+C+D=dI$, thus
\eqref{comrelbcdp} is satisfied.
Finally, since $BCD=(CB+X)D=CBD+XD=C(DB-X)+XD=CDB-CX+XD=CDB$, etc.,
one immediately deduces \eqref{comrelrot}.
\end{Remark}

We are ready to present our noncommutative extension of \eqref{7f6}:
\begin{Theorem}\label{nc7f6}
Let $A$, $B$, $C$ and $D$ be noncommutative parameters of some unit ring,
and assume that the relations \eqref{comrel} hold.
Further, let $n$ be a nonnegative integer. Then we have the following
summation for a noncommutative hypergeometric series of type I (and II).
\begin{multline}\label{nc7f6gl}
{}_7F_6\!\left\lceil\begin{matrix}\frac 12A+I,A,B,C,D,
2A-B-C-D+(n+1)I,\\\frac 12A,A+(n+1)I,A-C+I,A-B+I,
A-D+I,\end{matrix}\right.\\
\left.\begin{matrix}-nI\\B+C+D-A-nI\end{matrix}\,;I\right\rfloor\\=
\left\lceil\begin{matrix}A-C-D+I,A-B-D+I,A+I,A-B-C+I\\
A-C+I,A-B+I,A-D+I,A-B-C-D+I\end{matrix}\right\rfloor_n.
\end{multline}
\end{Theorem}

\begin{Remark}\label{remIaII}
The above summation holds for both types, I and II, of noncommutative
series and shifted factorials. More precisely, we could switch
the type I brackets to type II brackets on either side (or on both
sides) of \eqref{nc7f6gl} and the formula would be still valid.
This is a consequence of the conditions \eqref{comrel} (from which we
extracted, in particular, \eqref{comrelbcdd}).
Nevertheless, for brevity of display we write \eqref{nc7f6gl}
using type I brackets only.
\end{Remark}

\begin{proof}[Proof of Theorem~\ref{nc7f6}]
We prove \eqref{nc7f6gl} by induction on $n$.
For $n=0$ the formula is trivial.
Assume the formula is true up to a fixed $n$.
To prove it for $n+1$, use the elementary identity
\begin{multline*}
\left\lceil\begin{matrix}2A-B-C-D+(n+2)I,-(n+1)I\\
A+(n+2)I,B+C+D-A-(n+1)I\end{matrix}\right\rfloor_k\\=
\left\lceil\begin{matrix}2A-B-C-D+(n+1)I,-nI\\
A+(n+1)I,B+C+D-A-nI\end{matrix}\right\rfloor_k\\\times
\big[I-(2A-B-C-D+(n+1)I)^{-1}(2A-B-C-D+(2n+2)I)\\\times
(B+C+D-A-(n+1)I)^{-1}(B+C+D-A)\\\times
(A+(n+1+k)I)^{-1}(A+kI)((-n-1+k)I)^{-1}kI\big]
\end{multline*}
to obtain
\begin{multline*}
{}_7F_6\!\left\lceil\begin{matrix}
\frac 12A+I,A,B,C,D,
2A-B-C-D+(n+2)I,\\\frac 12A,A+(n+2)I,A-C+I,A-B+I,A-D+I,
\end{matrix}\right.\\
\left.\begin{matrix}
-(n+1)I\\B+C+D-A-(n+1)I\end{matrix}\,;I\right\rfloor\\=
\sum_{k=0}^n\left\lceil\begin{matrix}\frac 12A+I,A,B,C,D,
2A-B-C-D+(n+1)I,\\\frac 12A,A+(n+1)I,A-C+I,A-B+I,A-D+I,
I,\end{matrix}\right.\\
\left.\begin{matrix}-nI\\
B+C+D-A-nI\end{matrix}\right\rfloor_k\\\times
\big[I-(2A-B-C-D+(n+1)I)^{-1}(2A-B-C-D+(2n+2)I)\\\times
(B+C+D-A-(n+1)I)^{-1}(B+C+D-A)\\\times
(A+(n+1+k)I)^{-1}(A+kI)((-n-1+k)I)^{-1}kI\big]\\=
{}_7F_6\!\left\lceil\begin{matrix}\frac 12A+I,A,B,C,D,
2A-B-C-D+(n+1)I,\\
\frac 12A,A+(n+1)I,A-C+I,A-B+I,A-D+I,\end{matrix}\right.\\
\left.\begin{matrix}-nI\\
B+C+D-A-nI\end{matrix}\,;I\right\rfloor\\
-{}_7F_6\!\left\lceil\begin{matrix}\frac 12A+2I,A+2I,B+I,C+I,D+I,
2A-B-C-D+(n+2)I,\\\frac 12A+I,A+(n+3)I,A-C+2I,A-B+2I,
A-D+2I,\end{matrix}\right.\\
\left.\begin{matrix}-nI\\
B+C+D-A-(n-1)I\end{matrix}\,;I\right\rfloor\\\times
(B+C+D-A-(n+1)I)^{-1}(B+C+D-A)(A+(n+1)I)^{-1}(A+I)\\\times
(A+(n+2)I)^{-1}(A+2I)(B+C+D-A-nI)^{-1}(2A-B-C-D+(2+2n)I)\\\times
(A-C+I)^{-1}B(A-B+I)^{-1}C(A-D+I)^{-1}D\\=
\left\lceil\begin{matrix}A-C-D+I,A-B-D+I,A+I,A-B-C+I\\
A-C+I,A-B+I,A-D+I,A-B-C-D+I\end{matrix}\right\rfloor_n\\-
\left\lceil\begin{matrix}A-C-D+I,A-B-D+I,A+3I,A-B-C+I\\
A-C+2I,A-B+2I,A-D+2I,A-B-C-D\end{matrix}\right\rfloor_n\\\times
(B+C+D-A-(n+1)I)^{-1}(B+C+D-A)(A+(n+1)I)^{-1}(A+I)\\\times
(A+(n+2)I)^{-1}(A+2I)(B+C+D-A-nI)^{-1}(2A-B-C-D+(2+2n)I)\\\times
(A-C+I)^{-1}B(A-B+I)^{-1}C(A-D+I)^{-1}D,
\end{multline*}
the last equation due to the inductive hypothesis.
We are done after application of Lemma~\ref{lem76}.
\end{proof}

A natural question is if one can derive
a noncommutative terminating ${}_5F_4$ summation
(extending  \cite[Appendix~(III.13)]{S}) as a special case from
Theorem~\ref{nc7f6}. The answer is negative. If, say, $D=A-C+I$
then since $D$ and $C$ would commute, \eqref{comrelbcd}
would imply that all variables commute so one just obtains
the usual (commutative) terminating ${}_5F_4$ summation.

We now turn to the question of deriving a $Q$-analogue of
Theorem~\ref{nc7f6}.
The $q$-analogue of \eqref{7f6} is the following summation formula
(cf.\ \cite[Appendix~(II.22)]{GR}).
\begin{multline}\label{7f6q}
{}_8\phi_7\!\left[\begin{matrix}qa^{\frac 12},qa^{-\frac 12},a,
b,c,d,a^2q^{n+1}/bcd,q^{-n}\\
a^{\frac 12},a^{-\frac 12},aq/b,aq/c,aq/d,bcdq^{-n}/a,
aq^{n+1}\end{matrix};q,q\right]\\=
\frac{(aq;q)_n(aq/bc;q)_n(aq/bd;q)_n(aq/cd;q)_n}
{(aq/b;q)_n(aq/c;q)_n(aq/d;q)_n(aq/bcd;q)_n}.
\end{multline}

In order to derive a noncommutative extension of \eqref{7f6q},
one should at least be able to extend its $n=1$ special case,
which (with $a$ replaced by $a/q$) is
\begin{multline}\label{7f6q1}
1-\frac{(1-b)(1-c)(1-d)(1-a^2/bcd)}{(1-a/b)(1-a/c)(1-a/d)(1-bcd/a)}\\=
\frac{(1-a)(1-a/bc)(1-a/bd)(1-a/cd)}{(1-a/b)(1-a/c)(1-a/d)(1-a/bcd)}.
\end{multline}

A noncommutative extension of \eqref{7f6q1} is indeed available
when one assumes that $A$, $B$, $C$, and $D$ are parameters such that
$A$, $B$ and $D$ mutually commute while $C$ does not commute with
any of the other parameters. We then have
\begin{multline}\label{easy76Q}
I-(I-C^{-1}A)^{-1}(I-B)(I-AB^{-1})^{-1}(I-B^{-1}CB)\\\times
(I-DA^{-1}CB)^{-1}(I-D)(I-AD^{-1})^{-1}(I-AB^{-1}D^{-1}C^{-1}A)\\
=(I-C^{-1}A)^{-1}(I-AB^{-1}D^{-1})(I-AB^{-1})^{-1}(I-C^{-1}AD^{-1})\\\times
(I-B^{-1}C^{-1}AD^{-1})^{-1}(I-A)(I-AD^{-1})^{-1}(I-B^{-1}C^{-1}A).
\end{multline}
This is a consequence of \eqref{easy32Q} and the identity
\begin{multline}\label{easy32Q1}
I-D^{-1}(I-BA^{-1})^{-1}(I-B)(I-DA^{-1}CB)^{-1}D(I-A^{-1}CBDA^{-1})\\
=D^{-1}(I-BA^{-1})^{-1}(I-C^{-1}AD^{-1})(I-B^{-1}C^{-1}AD^{-1})^{-1}
(I-A^{-1})D,
\end{multline}
where $A$, $B$ and $D$ mutually commute while $C$ does not commute with
any of the other parameters.

Since our attempts to extend \eqref{easy76Q} to a noncommutative
extension of \eqref{7f6q} for general $n$ failed, we omit writing out the
technical details of the proofs of \eqref{easy32Q1} and \eqref{easy76Q}.
When trying to establish a $Q$-analogue of Theorem~\ref{nc7f6},
we had no problem with applying induction, when assuming that
$A$ and $DCB$ commute with all other parameters.
However, for the last step in the corresponding $Q$-analogue of
Lemma~\ref{lem76} to work out, we would need that
the three products $(I-B)(I-B^{-1}CB)(I-D)$,
$(I-AB^{-1}Q)(I-AC^{-1}Q)(I-AD^{-1}Q^2)$ and
$(I-AC^{-1}D^{-1}Q)(I-AD^{-1}B^{-1}Q)(I-AB^{-1}C^{-1}Q)$ mutually commute.
This means that additional conditions on the parameters
are needed. We were not able to specify any reasonable
additional conditions on $A$, $B$, $C$ and $D$ satisfying all of the
above but where $A$, $B$, $C$ and $D$ do not already  mutually commute.

The question remains whether there exists any more suitable
noncommutative extension of \eqref{7f6q1} than \eqref{easy76Q},
for the purpose of deriving a full $Q$-analogue of Theorem~\ref{nc7f6}.

\section{Summations for nonterminating series}\label{secnt}

The noncommutative (basic) hypergeometric series in
Sections~\ref{sec2f1}--\ref{sec7f6} terminate due to the occurrence
of the upper parameter $-nI$ (or $Q^{-n}$) where $n$ is a nonnegative
integer. Assuming that (the unit ring) $R$ is a Banach algebra with norm
$\|\cdot\|$, we may obtain summations for nonterminating noncommutative
(basic) hypergeometric series by (possibly substituting some variables
and) formally considering the term-wise limit as $n\to\infty$.
However, the validity of this procedure
which involves the interchange of limit and summation
would need to be justified case by case.
The main problem appears with convergence, especially for
noncommutative ``ordinary'' hypergeometric series,
for which precise conditions would need to be given
(as the argument in the summations of interest is often $I$).
Furthermore, in the limit a noncommutative generalized gamma
function (which would reduce to a product of quotients of
gamma functions in the case of commuting parameters)
would be needed which we did not (yet) define.
This problem does not appear in the basic $Q$-case where
infinite $Q$-shifted factorials make perfectly sense assuming
$\|Q\|<1$. Wherever we consider infinite $Q$-shifted factorials
we shall implicitly assume such $Q$.

The following two identities are obtained from Theorem~\ref{nc3f2Q}
by taking the limit $n\to\infty$ on each side, termwise on the
left-hand sides. These are conjectured noncommutative $Q$-Gau{\ss}
summations, extending \cite[Appendix~(II.8)]{GR}.

\begin{Conjecture}\label{ncnt2f1Q}
Let $A$, $B$ and $C$ be noncommutative parameters of some Banach algebra,
and suppose that $Q$ commutes each with $A$, $B$ and $C$. Further, assume
that the product $BC^{-1}A$ commutes each with $A$, $B$ and $C$.
Moreover, assume that $\|A^{-1}CB^{-1}\|<1$.
Then we have the following summation for a noncommutative basic
hypergeometric series of type I.
\begin{equation}\label{ncnt2f1QIgl}
{}_2\phi_1\!\left\lceil\begin{matrix}A,B\\
C\end{matrix}\,;Q,A^{-1}CB^{-1}\right\rfloor=
\left\lceil\begin{matrix}CB^{-1},A^{-1}C\\
C,A^{-1}CB^{-1}\end{matrix};Q,I\right\rfloor_\infty.
\end{equation}
Further, we have the following
summation for a noncommutative basic hypergeometric series of type II.
\begin{align}\label{ncnt2f1QIIgl}\notag
{}_2\phi_1\!\left\lfloor\begin{matrix}A,B\\
C\end{matrix}\,;Q,A^{-1}CB^{-1}\right\rceil={}&
\left\lfloor\begin{matrix}CB^{-1}Q,A^{-1}CQ\\
C,A^{-1}CB^{-1}\end{matrix};Q,I\right\rceil_\infty\\&\times
(I-CB^{-1})(I-A^{-1}C).
\end{align}
\end{Conjecture}

Next we give two noncommutative
extensions of the nonterminating $q$-binomial theorem (cf.\
\cite[II.3]{GR}).

\begin{Theorem}\label{nc1f0Q}
Let $A$ and $C$ be noncommutative parameters of some Banach algebra,
and suppose that $Q$ commutes both with $A$ and $C$. Further,
assume that $\|Z\|<1$. Then we have the following
summation for a noncommutative hypergeometric series of type I.
\begin{equation}\label{nc1f0QIgl}
{}_1\phi_0\!\left\lceil\begin{matrix}A\\
-\end{matrix};Q,Z\right\rfloor=
\left\lfloor\begin{matrix}AZ\\Z\end{matrix};Q,I\right\rceil_\infty.
\end{equation}
Further, we we have the following
summation for a noncommutative hypergeometric series of type II.
\begin{equation}\label{nc1f0QIIgl}
{}_1\phi_0\!\left\lfloor\begin{matrix}A\\
-\end{matrix};Q,Z\right\rceil=
\left\lceil\begin{matrix}AZQ^{-1}\\Z\end{matrix};Q,I\right\rfloor_\infty
(I-AZQ^{-1})^{-1}.
\end{equation}
\end{Theorem}
Note that with the notation defined in Subsection~\ref{subsecrp},
\eqref{nc1f0QIIgl} can be written more compactly as
\begin{equation}
{}_1\phi_0\!\left\lfloor\begin{matrix}A\\
-\end{matrix};Q,Z\right\rceil=
{}^{\text{\raisebox{1em}{$\sim$}}}\!
\left\lfloor\begin{matrix}AZ\\Z\end{matrix};Q,I\right\rceil_\infty.
\end{equation}

Identities \eqref{nc1f0QIgl} and \eqref{nc1f0QIIgl} can be obtained
from Theorem~\ref{nc2f1Q} by performing the substitution
$C\mapsto Q^{1-n}Z^{-1}$, taking the limit $n\to\infty$ on each side
of the respective identities, and multiplying each side
by $Z^{-1}$ from the left and by $Z$ from the right.
Since we did not justify taking termwise limits, we provide
an independent proof.
\begin{proof}[Proof of Theorem~\ref{nc1f0Q}]
We prove \eqref{nc1f0QIgl}, leaving the proof of \eqref{nc1f0QIIgl}
to the reader.

Let $f(A,Z)$ denote the series on the left-hand side of \eqref{nc1f0QIgl}.
We make use of the two simple identities
\begin{subequations}
\begin{align}\label{fq1}
Z={}&AZQ^k+(I-AQ^k)Z,\\\label{fq2}
I={}&Q^k+(I-Q^k),
\end{align}
\end{subequations}
to obtain two functional equations for $f$.
First, \eqref{fq1} gives
\begin{equation}\label{fq11}
Zf(A,Z)=AZf(A,ZQ)+f(AQ,Z)(I-A)Z,
\end{equation}
while \eqref{fq2} gives
\begin{equation}\label{fq21}
f(A,Z)=f(A,ZQ)+f(AQ,Z)(I-A)Z.
\end{equation}
Combining \eqref{fq11} and \eqref{fq21}, one immediately has
\begin{equation*}
f(A,Z)=f(A,ZQ)+Zf(A,Z)-AZf(A,ZQ),
\end{equation*}
or equivalently
\begin{equation*}
(I-Z)f(A,Z)=(I-AZ)f(A,ZQ),
\end{equation*}
thus
\begin{equation}\label{fq0}
f(A,Z)=(I-Z)^{-1}(I-AZ)f(A,ZQ).
\end{equation}
Iteration of \eqref{fq0} gives the result since $f(A,O)=I$ by definition
of $f$.
\end{proof}

\section{More identities}\label{secmore}

\subsection{Telescoping}

By iterating any of the simple identities \eqref{easy32},
\eqref{easy32Q}, \eqref{easy76} or \eqref{easy76Q}, indefinite
summations involving general parameters can be derived.
For example, using \eqref{easy32Q} one has
\begin{multline}
\sum_{k=0}^n\prod_{j=0}^{k-1}\left[
(I-C_j)^{-1}(I-A_j)(I-B_jC_j^{-1}A_j)^{-1}(I-B_j)\right]\\\times
(I-C_k)^{-1}(I-C_kB_k^{-1})(I-A_k^{-1}C_kB_k^{-1})^{-1}(I-A_k^{-1}C_k)\\
=\sum_{k=0}^n\left[\prod_{j=0}^{k-1}\left[
(I-C_j)^{-1}(I-A_j)(I-B_jC_j^{-1}A_j)^{-1}(I-B_j)\right]\right.\\
-\left.\prod_{j=0}^k\left[
(I-C_j)^{-1}(I-A_j)(I-B_jC_j^{-1}A_j)^{-1}(I-B_j)\right]\right]\\
=I-\prod_{j=0}^n\left[
(I-C_j)^{-1}(I-A_j)(I-B_jC_j^{-1}A_j)^{-1}(I-B_j)\right],
\end{multline}
for arbitrary noncommutative parameters $A_j$, $B_j$, $C_j$,
where $0\le j\le n$.

One can now obtain a basic (or even multi-basic) identity in a natural way
by setting $A_j=AQ^j$ etc.\ and similarly for $B_j$ and $C_j$.
We leave the details to the reader.

\subsection{Reversing all products}\label{subsecrp}

For all identities given in this paper, one may obtain new ones
by simply reversing all the products (of elements of the unit ring $R$)
simultaneously on each side of the respective identities.
This is clearly an involution.
(For square matrices it would amount to transposition of matrices.)

For instance, with the new definitions
\begin{equation}\label{defncpochr}
{}^{\text{\raisebox{1em}{$\sim$}}}\!\left\lceil\begin{matrix}
A_1,A_2,\dots,A_r\\C_1,C_2,\dots,C_r
\end{matrix};Z\right\rfloor_k:=
\prod_{j=1}^k\left(Z\prod_{i=1}^r(A_i+(j-1)I)
(C_i+(j-1)I)^{-1}\right),
\end{equation}
and
\begin{equation}\label{defnchypr}
{}_{r+1}F_r\!\!^{\text{\raisebox{1em}{$\sim$}}}\!\left\lceil\begin{matrix}
A_1,A_2,\dots,A_{r+1}\\C_1,C_2,\dots,C_r
\end{matrix};Z\right\rfloor:=
\sum_{k\ge 0}
{}^{\text{\raisebox{1em}{$\sim$}}}\!\left\lceil\begin{matrix}
A_1,A_2,\dots,A_{r+1}\\C_1,C_2,\dots,C_r,I,
\end{matrix};Z\right\rfloor_k
\end{equation}
for reversed (or ``transposed'') versions of
generalized noncommutative shifted factorials and
noncommutative hypergeometric series of type I
(compare with \eqref{defncpochI} and \eqref{defnchypI};
similar definitions can be made for type II and in the basic cases),
and writing
\begin{equation}\label{defncpochar}
{}^{\text{\raisebox{1em}{$\sim$}}}\!\left\lceil\begin{matrix}
A_1,A_2,\dots,A_r\\C_1,C_2,\dots,C_r
\end{matrix};I\right\rfloor_k=
{}^{\text{\raisebox{1em}{$\sim$}}}\!\left\lceil\begin{matrix}
A_1,A_2,\dots,A_r\\C_1,C_2,\dots,C_r
\end{matrix}\right\rfloor_k
\end{equation}
for brevity,
we have the following noncommutative Chu--Vandermonde summation:
\begin{Theorem}\label{nc2f1Ir}
Let $A$, $C$ be noncommutative parameters of some unit ring. Then
\begin{equation}\label{nc2f1Irgl}
{}_2F_1\!\!^{\text{\raisebox{1em}{$\sim$}}}\!
\left\lceil\begin{matrix}A,-nI\\C\end{matrix};I\right\rfloor=
{}^{\text{\raisebox{1em}{$\sim$}}}\!
\left\lceil\begin{matrix}C-A\\C\end{matrix}\right\rfloor_n.
\end{equation}
\end{Theorem}
This is a simple consequence of Theorem~\ref{nc2f1}, obtained by
reversing all products. Similarly, all the other identities
involving noncommutative parameters appearing in this paper
have reversed versions. We do not write them out explicitly.

\end{document}